\documentclass[11pt,a4paper,reqno]{amsart}
\usepackage[T1]{fontenc}
\usepackage[utf8]{inputenc}
\usepackage{lmodern,amsmath,amssymb,amsthm,microtype,tikz,needspace,placeins}
\usepackage[left=1in,right=1in,top=1in,bottom=1in]{geometry}
\usepackage[unicode,colorlinks=true,linkcolor=black,citecolor=blue!45!black,urlcolor=blue!45!black]{hyperref}
\hypersetup{pdfauthor={Yutian Li},pdftitle={A proof of the polygonal Faber–Krahn inequality},bookmarksdepth=2}
\newtheorem{theorem}{Theorem}[section]
\newtheorem{proposition}[theorem]{Proposition}
\newtheorem{lemma}[theorem]{Lemma}
\newtheorem{corollary}[theorem]{Corollary}
\theoremstyle{remark}
\newtheorem*{remark}{Remark}
\allowdisplaybreaks[2]
\title[A proof of the polygonal Faber--Krahn inequality]{A proof of the polygonal Faber--Krahn inequality}
\author{Yutian Li}
\address{School of Mathematics and Statistics, Nanfang College, Guangzhou, 510970 Guangdong, China.}
\email{yutianlee@gmail.com}
\date{20 November 2026}
\hypersetup{pdfkeywords={polygonal Faber–Krahn inequality, Pólya–Szegő conjecture, Dirichlet eigenvalues, shape optimization, mixed boundary conditions, spectral zeta functions, spectral partitions}}
\begin{document}
\begin{abstract}
We prove that, for every integer \(n\ge3\), the regular \(n\)-gon uniquely minimizes the first Dirichlet eigenvalue, up to congruence, among simple polygons of prescribed area with at most \(n\) sides. We also obtain the optimal linear correction when fewer than \(n\) sides are present. The proof combines a partition into convex cells with a mixed boundary comparison. Among convex polygons of fixed area and fixed sum of base angles less than \(\pi\), the isosceles triangle minimizes the first eigenvalue with Dirichlet conditions on the base and Neumann conditions elsewhere. A comparison of Fourier transforms of conformal area densities proves this inequality and the corresponding upper bounds for inverse spectral moments of integer order at least two. We prove strict log-convexity of an interpolation of the eigenvalues of regular polygons of unit area for real parameters at least four. We then derive sharp hexagonal bounds for periodic convex polygonal tilings and asymptotic bounds for convex polygonal partitions of a fixed convex polygon, with an explicit correction for boundary junctions.
\end{abstract}
\maketitle

\section{Introduction}\label{sec:1}

Among planar domains of prescribed area, the disk minimizes the first Dirichlet eigenvalue, as Faber \cite{ref26} and Krahn \cite{ref27} proved. For polygons of prescribed area with \(n\) sides, Pólya and Szegő conjectured in 1951 that the minimizer is the regular \(n\)-gon \cite[p. 158]{ref23}. Steiner symmetrization gives the equilateral triangle and the square in the first two cases \cite{ref22,ref23}. For a polygon with five or more sides, symmetrization can increase the number of sides \cite{ref3}.

For a bounded open domain \(\Omega\subset\mathbb R^2\), let \(|\Omega|\) denote its area and \(\lambda_j(\Omega)\) the \(j\)th Dirichlet eigenvalue of \(-\Delta\), counted in nondecreasing order with multiplicity. By the Rayleigh principle,
\[
 \lambda_1(\Omega)=\inf_{0\ne u\in H_0^1(\Omega)}
     \frac{\int_\Omega|\nabla u|^2\,dx}{\int_\Omega|u|^2\,dx}.
\]
The product \(\Phi(\Omega)=|\Omega|\lambda_1(\Omega)\) is invariant under similarities. A simple polygon \(P\) is the bounded domain enclosed by a simple polygonal curve. A side is a maximal straight segment of \(\partial P\), and interior angles may exceed \(\pi\). For each integer \(n\ge3\), let \(R_n\) be a regular \(n\)-gon and put \(c_n=\Phi(R_n)\).

\begin{theorem}[Polygonal Faber–Krahn inequality]\label{th:1.1}
Every simple polygon \(P\) with \(k\) sides satisfies \(\Phi(P)\ge c_k\), with equality exactly for a regular polygon. For \(n\ge4\), let \(\Delta_n=c_{n-1}-c_n\). If \(3\le k\le n\), then
\[
 \Phi(P)\ge c_n+(n-k)\Delta_n,\qquad \Delta_n>0. \tag{1.1}\label{eq:1.1}
\]
The coefficient \(\Delta_n\) is optimal. Equality in \eqref{eq:1.1} holds precisely for polygons similar to \(R_n\) or \(R_{n-1}\). Thus among polygons of prescribed area with at most \(n\) sides, the regular \(n\)-gon uniquely minimizes the eigenvalue up to congruence.
\end{theorem}

Faifman and Polterovich proved minimality on the area-preserving affine orbit of each regular polygon \cite{ref14}, and Indrei obtained minimality for an explicit class of competitors when the side count is sufficiently large \cite{ref18}. For fixed \(n\ge5\), Bogosel and Bucur showed that either certified computation can decide the conjecture or a nonregular polygon in the closure of the admissible class has the same value as the regular polygon \cite[Theorem 7.6]{ref3}. They subsequently established local minimality of the regular pentagon and hexagon by validated computation \cite{ref4}. Cheng, Gui, Hu and Li treated all convex polygons with sufficiently many sides \cite{ref11}. Wang \cite{ref25} proves the inequality for all \(n\ge5\), as well as results on torsional rigidity, Sobolev quotients and stability.

For a polygon with a designated side, called the base, write \(\lambda_j^{\mathrm{DN}}\) for the \(j\)th eigenvalue with Dirichlet conditions on the base and Neumann conditions on the remaining sides, counted with multiplicity in nondecreasing order. For \(0<\alpha<\pi\), let \(I_\alpha\) be the isosceles triangle with altitude one and apex angle \(\alpha\), with the opposite side as its base, and put
\[
 F(\alpha)=|I_\alpha|\lambda_1^{\mathrm{DN}}(I_\alpha).
\] Joining the center of \(R_n\) to its vertices divides it into \(n\) congruent triangles. The positive first eigenfunction, or ground state, of \(R_n\) has zero normal derivative on the radial segments by reflection symmetry. Its restriction to each triangle is positive and hence gives the first mixed eigenvalue there. By scaling,
\[
 c_n=nF(2\pi/n).
\]

The proof uses three ingredients. For a convex polygon \(Q\) with a designated base, let \(\gamma\) be \(\pi\) minus the sum of the two interior angles adjacent to the base. For a triangle, \(\gamma\) is its apex angle. Proposition~\ref{th:2.1} gives \(|Q|\lambda_1^{\mathrm{DN}}(Q)\ge F(\gamma)\) whenever \(\gamma>0\). Proposition~\ref{th:4.1} shows that, for \(n\ge4\), the tangent to \(F\) at \(2\pi/n\) supports its graph for every \(0\le\gamma<\pi\). Finally, Lemma~\ref{th:5.1} partitions an arbitrary simple polygon into convex cells with bases on the original boundary. In each cell, \(\gamma\) exceeds the horizontal intercept of the tangent by an amount proportional to the cell area. These three statements give the same lower bound for the first mixed eigenvalue on every cell. Restricting a Dirichlet test function to the cells and summing the energy inequalities proves the polygonal bound.

Triangle covers also enter the proof by Solynin and Zalgaller that the regular polygon minimizes logarithmic capacity at fixed area and side count \cite{ref24}. Our construction uses the partitions of Bogosel, Bucur and Fragalà, who relate cell areas to angles in their comparison for torsional rigidity \cite{ref5,ref6}. Wang uses the Kohler–Jobin inequality to pass from a torsion deficit to the eigenvalue bound for nonconvex polygons \cite{ref25}. His rearrangement on convex cells gives the mixed first-eigenvalue comparison \cite[Lemma 2.1]{ref25}. The Fourier comparison in Section~\ref{sec:2} also gives an upper bound for \(\sum_{j\ge1}(|Q|\lambda_j^{\mathrm{DN}}(Q))^{-r}\) at every integer \(r\ge2\), with equality only for the corresponding isosceles triangle (Corollary~\ref{th:6.7}). For Fourier rearrangement in Euclidean space, see \cite{ref19}.

We show that \(x\mapsto\log[xF(2\pi/x)]\) has positive second derivative for \(x\ge4\). Together with the comparison between the triangle and the square, this gives the monotonicity of \(c_j/c_{j+1}\) proved by Dahne, Gómez-Serrano and Pech-Alberich \cite{ref12}, following the monotonicity results of Antunes and Freitas \cite{ref1}. For a partition, Euler's formula bounds the mean number of sides. Applying a supporting line at six gives a lower bound for the geometric mean of the cell eigenvalues that accounts for nonhexagonal cells, irregular shapes, and unequal areas. For a partition of a bounded convex polygon, there is also a term for junctions on the outer boundary (Theorems~\ref{th:6.3} and~\ref{th:6.4}).

\subsection{Proof of the main theorem}

\noindent\textit{Proof of Theorem~\ref{th:1.1}.}\enspace For \(n=3\), the result is the classical minimizing property of the equilateral triangle \cite{ref22,ref23}. Appendix~\ref{sec:A} gives a proof using the boundary-moment argument of Fragalà and Velichkov \cite{ref29}.

In coordinates the vertices of \(I_\alpha\) are \((0,0),(1,\pm t)\), where \(t=\tan(\alpha/2)\). The trial function \(1-x\) gives \(0<F(\alpha)\le6t\), so \(F\) extends continuously to \(F(0)=0\).

Now fix \(n\ge4\) and let \(P\) have \(3\le k\le n\) sides. Put \(a=2\pi/n\) and \(s=F'(a)\). By Proposition~\ref{th:4.1}, \(F(a)-as<F(0)=0\), so \(s>F(a)/a>0\). The tangent has value \(s(\gamma-b)\) at angle \(\gamma\), where \(b=a-F(a)/s\) is its horizontal intercept. In particular \(0<b<a\). To bound \(b\), map \(I_{\pi/2}\) onto \(I_\alpha\) by \((x,Y)\mapsto(x,tY)\). This gives
\[
 \frac{F(\alpha)}{t}
 =\inf_{\substack{0\ne v\in H^1(I_{\pi/2})\\ v|_{x=1}=0}}
 \frac{\int_{I_{\pi/2}}(|v_x|^2+t^{-2}|v_Y|^2)}
      {\int_{I_{\pi/2}}|v|^2}.
\]
The admissible functions are independent of \(t\), and each quotient is nonincreasing in \(t\). Since \((d/d\alpha)\log t=1/\sin\alpha\), we obtain
\[
 \sin a\,F'(a)\le F(a),\qquad b\le a-\sin a .
\]
Since \(a\le\pi/2\), concavity of sine gives \(\sin a\ge2a/\pi>a/2\), so \((k-1)b<kb<na/2=\pi\). Put \(c=(2\pi-kb)/|P|\).

By Lemma~\ref{th:5.1}, \(P\) is partitioned, up to the cell boundaries, into convex polygons \(V_j\). Each has a base that is a whole side of \(V_j\) and lies on \(\partial P\). Writing \(\gamma_j\) for \(\pi\) minus the sum of its two base angles, we have \(\gamma_j\ge b+c|V_j|>0\). Thus
\[
 |V_j|\lambda_1^{\mathrm{DN}}(V_j)\ge F(\gamma_j)
 \ge s(\gamma_j-b)\ge sc|V_j|. \tag{1.2}\label{eq:1.2}
\]
Every \(u\in H_0^1(P)\) restricts to a function with zero trace on the base of each \(V_j\). Summing the resulting energy inequalities gives \(\lambda_1(P)\ge sc\), hence
\[
 \Phi(P)\ge s(2\pi-kb)=nF(a)+(n-k)sb. \tag{1.3}\label{eq:1.3}
\]
Write \(d_n=sb=aF'(a)-F(a)>0\). In particular, \(\Phi(P)\ge c_n\) when \(k=n\ge4\).

\emph{Equality.} If \(\Phi(P)=c_n\), then \eqref{eq:1.3} forces \(k=n\). The ground state has positive mass on each cell, so equality in the sum forces equality in \eqref{eq:1.2} on every cell. The restriction of the ground state minimizes the mixed Rayleigh quotient there. The tangent touches the graph only at \(a\), so \(\gamma_j=a\), and equality in Proposition~\ref{th:2.1} makes \(V_j\) isosceles. Each triangle has area \((a-b)/c=|P|/n\), so there are exactly \(n\) congruent triangles.

The bases cover \(\partial P\) apart from finitely many endpoints. Indeed, if a Neumann leg lay along a boundary segment, the positive mixed ground state would vanish there by the Dirichlet condition on \(P\). Even reflection across the leg and the strong maximum principle exclude this. Every side of \(P\) therefore contains a base, and each base lies on a single side. There are \(n\) of each, so the bases are whole sides and all sides have equal length. At each vertex the two incident triangles occupy disjoint sectors with angles totaling \(\pi-a\). These lower bounds sum to \((n-2)\pi\), the angle sum of \(P\), so every angle is \(\pi-a\). Thus \(P\) is regular. Conversely, the regular polygon attains the bound.

\emph{The optimal coefficient.} Compare \(F(2\pi/i)\) with its lower bound from the tangent at \(2\pi/j\) and multiply by \(i\). Proposition~\ref{th:4.1} gives, for \(j\ge4\),
\[
 c_i>c_j+(j-i)d_j\qquad(i\ge3,\ i\ne j).
\]
Taking \(i=j-1\) and \(i=j+1\) gives \(\Delta_j>d_j>\Delta_{j+1}\). In particular every \(\Delta_j\) with \(j\ge4\) is positive. The fixed-side comparison and \(c_k-c_n=\sum_{j=k+1}^n\Delta_j\) now prove \eqref{eq:1.1}. If \(k\le n-2\), the sum is strictly greater than \((n-k)\Delta_n\). For \(k=n\) and \(k=n-1\), equality holds precisely for the regular polygon with that many sides. The case \(k=n-1\) shows that the coefficient is optimal.\qed

\section{Fourier comparison for a mixed boundary problem}\label{sec:2}

We write the mixed Rayleigh quotients on a fixed strip and compare the resulting area densities by Fourier transform. The exponential map \((x,t)\mapsto e^{x+it}\) sends the strip \(S=\mathbb R\times(0,\pi)\) onto the upper half-plane. We extend densities on \(S\) evenly across \(t=0\) and then periodically to the cylinder \(\mathcal C=\mathbb R\times(\mathbb R/2\pi\mathbb Z)\). For an integrable density \(p\) on \(\mathcal C\), we use the Fourier transform
\[
 \widehat p(\xi,k)=\int_{\mathcal C} e^{i(\xi x+kt)}p(x,t)\,d\mu,\qquad
 d\mu=\frac{dx\,dt}{2\pi},\qquad (\xi,k)\in\mathbb R\times\mathbb Z .
\]

\begin{proposition}[Comparison with an isosceles triangle]\label{th:2.1}
Let \(Q\) be a bounded convex polygon with one side designated as the base. Use the space of functions in \(H^1(Q)\) with zero trace on the base, so the other sides have the natural Neumann condition. Suppose the angles at the base endpoints are \(\pi a\) and \(\pi b\), with \(a,b>0\) and \(C=a+b<1\). Then
\[
 |Q|\lambda_1^{\mathrm{DN}}(Q)\ge F(\pi(1-C)), \tag{2.1}\label{eq:2.1}
\]
Equality holds exactly when \(Q\) is similar to \(I_{\pi(1-C)}\) with the bases corresponding. Let \(p\) and \(p_*\) be the normalized pullbacks of area on \(Q\) and \(I_{\pi(1-C)}\), extended evenly to \(\mathcal C\), by the conformal maps defined in Section~\ref{sec:strip-map}. They are probability densities with respect to \(d\mu\). For the triangle, the apex is the image of \(1\) in the upper half-plane coordinates. Then \(|\widehat p(\xi,k)|\le\widehat p_*(\xi,k)\) at every frequency. The inequality is strict away from \((0,0)\) unless \(Q\) is similar to \(I_{\pi(1-C)}\) with corresponding base.
\end{proposition}

We first compare \(Q\) with the triangle having the same base angles. The density for \(Q\) is a convolution of the triangle density with a probability measure. We then compare the Fourier transforms of the two triangle densities. The following lemma passes from a bound on these transforms to an eigenvalue inequality.

\Needspace{8\baselineskip}
\begin{lemma}[Fourier comparison on the strip]\label{th:2.2}
Let \(p\) and \(q\) be nonzero nonnegative integrable densities on \(\mathcal C\), even in the angular variable. On the space of \(u\in H^1(S)\) with zero trace on \(t=\pi\), set
\[
 E(u)=\int_S|\nabla u|^2\,dx\,dt,\qquad
 M_p(u)=\int_Sp|u|^2\,\frac{dx\,dt}{\pi}.
\]
Define \(M_q\) by the same formula with \(q\) in place of \(p\). Assume that both mass forms are continuous in the norm \(E^{1/2}\). On this space define
\[
 \Lambda(p)=\inf_{M_p(u)>0}\frac{E(u)}{M_p(u)}.
\]
Define \(\Lambda(q)\) in the same way. If \(\widehat q(\xi,k)\ge|\widehat p(\xi,k)|\) at every frequency, then \(\Lambda(p)\ge\Lambda(q)\). If also the infimum for \(\Lambda(p)\) is attained and
\[
 \widehat q(\xi,0)>|\widehat p(\xi,0)|\qquad(\xi\ne0),
\]
then the eigenvalue inequality is strict.
\end{lemma}

\noindent\textit{Proof.}\enspace Expand \(u\) in the basis \(\sqrt2\cos((m+\tfrac12)t)\), \(m=0,1,\ldots\), orthonormal for \(dt/\pi\). Let \(U_m\) be the unitary Fourier transform in \(x\) of the coefficient \(u_m(x)\), with the convention \(U_m(\xi)=(2\pi)^{-1/2}\int_{\mathbb R}e^{-i\xi x}u_m(x)\,dx\) on integrable functions. By Parseval's identity the energy is diagonal, with multiplier \(\xi^2+(m+\tfrac12)^2\) and measure \(\pi\,d\xi\) on each mode. For angular indices \(m,n\ge0\) and real frequencies \(\xi,\eta\), multiplying the two cosines gives the coefficient of \(U_m(\xi)\overline{U_n(\eta)}\) in \(M_p(u)\) as
\[
 \frac{\widehat p(\xi-\eta,m-n)+\widehat p(\xi-\eta,m+n+1)}{2\pi}.
\]
Replace the coefficients \(U_m\) by \(|U_m|\) to obtain \(v\). This preserves \(E\) and the zero trace on \(t=\pi\). The triangle inequality and the bound on \(\widehat p\) give \(M_p(u)\le M_q(v)\) for finitely many angular modes with bounded horizontal frequency support. The same inequality holds for all admissible functions by energy convergence, since taking moduli is continuous in the energy norm and the mass forms are continuous. Taking infima gives \(\Lambda(p)\ge\Lambda(q)\).

For strictness, choose an index \(m\ge0\) for which \(U_m\not\equiv0\), and a bounded measurable set \(B\) of positive measure on which \(|U_m|\) is bounded below by a positive constant. Up to the factor \(1/(2\pi)\), the term with both angular indices equal to \(m\) contributes
\[
 \iint_{B\times B}
 [\widehat q(\xi-\eta,0)-|\widehat p(\xi-\eta,0)|]\,
 |U_m(\xi)U_m(\eta)|\,d\xi\,d\eta>0.
\]
The integrand is positive off the diagonal, which has measure zero. After taking real parts and applying the triangle inequality, every remaining contribution to \(M_q(v)-M_p(u)\) is nonnegative. Thus this strict bound holds for every truncation containing \(B\) in mode \(m\), and hence in the limit. Applying the bound to a minimizer for \(\Lambda(p)\) gives \(\Lambda(q)<\Lambda(p)\).\qed

\subsection{The conformal map from the strip}\label{sec:strip-map}

\noindent\textit{Proof of Proposition~\ref{th:2.1}.}\enspace Map the upper half-plane conformally onto \(Q\), taking \(0,\infty\) to the endpoints with angles \(\pi a,\pi b\), respectively, and the negative real axis to the base. Let \(r_1,\ldots,r_N>0\) be the preimages of the remaining vertices, and let \(\pi\delta_j\) be the exterior angle at the image of \(r_j\). The sides are maximal straight segments, so \(\delta_j>0\), and the angle sum gives \(\sum_j\delta_j=C\). By the Schwarz–Christoffel formula \cite{ref13}, the derivative is proportional to \(z^{a-1}\prod_j(z-r_j)^{-\delta_j}\). The exponents follow by straightening each corner with a power map. After division by this product of powers, the derivative has constant boundary argument modulo \(\pi\). Reflection and the behavior at infinity extend the quotient holomorphically to the sphere, so it is constant.

For \((x,t)\in S\) set \(z=e^{x+it}\). The Jacobian of the composite map from \(S\) to \(Q\) is proportional to \(|z|^{2a}\prod_j|z-r_j|^{-2\delta_j}\). Multiply the Jacobian by \(\pi/|Q|\) and extend it evenly to \(\mathcal C\) to obtain the probability density \(p\) for \(d\mu\). Conformal invariance of the Dirichlet integral gives
\[
 |Q|\lambda_1^{\mathrm{DN}}(Q)=\inf_{\substack{0\ne u\in H^1(S)\\ u(\cdot,\pi)=0}}
 \frac{\int_S|\nabla u|^2\,dx\,dt}{\int_S p|u|^2\,dx\,dt/\pi}.
 \tag{2.2}\label{eq:2.2}
\]
The Poincaré inequality in \(t\) identifies the finite-energy functions of zero trace on \(t=\pi\) with those in the infimum in \eqref{eq:2.2}, together with zero. The density decays exponentially at both infinite ends. Near the preimage of the vertex corresponding to \(r_j\), it has a singularity of order \(-2\delta_j>-2\). Thus \(p\in L^{1+\epsilon}\) for some \(\epsilon>0\), and Hölder's inequality and Sobolev embedding make the mass continuous in the energy norm. Logarithmic cutoffs at the base endpoints have energy tending to zero, so the endpoints impose no further trace conditions. The compact Sobolev embedding on the bounded polygon gives a minimizer.

\subsection{Convolution of the area densities}

Write \(p_{a,b}\) for the density of the triangle with base angles \(\pi a,\pi b\), whose apex is the image of \(1\). Writing \(*\) for convolution on \(\mathcal C\), we show that for some probability measure \(\nu\),
\[
                           p=p_{a,b}*\nu. \tag{2.3}\label{eq:2.3}
\]
Let \(\Gamma\) denote Euler's Gamma function. Apply its integral representation in the form
\[
 \Gamma(\delta_j)|z-r_j|^{-2\delta_j}
 =\int_0^\infty s_j^{\delta_j-1}e^{-s_j|z-r_j|^2}\,ds_j.
\] Change variables to the sum \(\sum_j s_j\) and the proportions \(t_j=s_j/\sum_\ell s_\ell\). Integrating out the sum expresses \(\prod_j|z-r_j|^{-2\delta_j}\) as the average of \((\sum_j t_j|z-r_j|^2)^{-C}\) on the simplex \(\{t_j\ge0,\ \sum_j t_j=1\}\). The probability density in the first \(N-1\) coordinates is \(\Gamma(C)\prod_j t_j^{\delta_j-1}/\prod_j\Gamma(\delta_j)\). If \(N=1\), use the point mass at \(t_1=1\). Put \(w_0=\sum_j t_jr_j\) and \(h=(\sum_j t_j|r_j-w_0|^2)^{1/2}\). Then \(\sum_j t_j|z-r_j|^2=|z-w_0|^2+h^2\).

For \(h>0\), an integral of translates of \(|z|^{-2C}\) gives
\[
 (|z-w_0|^2+h^2)^{-C}
 =\frac{1-C}{\pi}\int_{\mathbb C}|z-w|^{-2C}
       \frac{h^{\,2-2C}}{(|w-w_0|^2+h^2)^{2-C}}\,dA(w).
 \tag{2.4}\label{eq:2.4}
\]
Here \(dA\) denotes planar area. By translation and scaling, it suffices to take \(w_0=0\) and \(h=1\). Stereographic projection sends uniform probability on the unit sphere to \(d\omega(w)=dA(w)/[\pi(1+|w|^2)^2]\). Half the chord length between the points corresponding to \(z\) and \(w\) is \(d(z,w)=|z-w|/\sqrt{(1+|z|^2)(1+|w|^2)}\). By rotational invariance, \(\int d(z,w)^{-2C}\,d\omega(w)\) is independent of \(z\). At \(z=0\), polar integration shows that \(d(0,w)^2=|w|^2/(1+|w|^2)\) is uniform on \([0,1]\), so the integral is \(1/(1-C)\). This gives \eqref{eq:2.4}. For \(h=0\), use a point mass at \(w_0\) as the measure of translates.

For \(w\ne0\), the factor \(|z|^{2a}|z-w|^{-2C}\) is \(|w|^{-2b}\) times the triangle's density before normalization, translated on \(\mathcal C\) by \((\log|w|,\arg w)\). Tonelli's theorem and finiteness of the area allow us to normalize the weights of these translates, giving \eqref{eq:2.3}. If \(N>1\), then \(h>0\) in the interior of the simplex, and the measure of translates in \eqref{eq:2.4} has a positive density throughout the plane. Consequently \(\nu\) gives positive mass to every nonempty open set on \(\mathcal C\). Its Fourier transform is the integral of \(e^{i(\xi x+kt)}\) against \(\nu\). At a nonzero frequency this function takes different values on the support of \(\nu\), so equality in the triangle inequality is impossible and \(|\widehat\nu(\xi,k)|<1\).

\subsection{The transform of a triangle}

Put \(q=1-C\), \(u=(|k|+i\xi)/2\) and \(p_*=p_{C/2,C/2}\). Define \(r_u(A)\) for \(A>0\) by the Gamma ratio below. Then
\[
 r_u(A)=\frac{\Gamma(A+u)\Gamma(A+q)}{\Gamma(A)\Gamma(A+q+u)},\qquad
 \widehat p_{a,b}(\xi,k)=r_u(a)r_{\bar u}(b).
 \tag{2.5}\label{eq:2.5}
\]
The triangle density is proportional to \(e^{(a-b)x}(\cosh x-\cos t)^{-C}\). For real \(A>B\), the beta integral gives
\[
 (A-B)^{-C}=\frac{\sin(\pi C)}{\pi}
       \int_A^\infty\frac{(v-A)^{-C}}{v-B}\,dv.
\]
Take \(A=\cosh x,\ B=\cos t,\ v=\cosh w\). The angular factor
\(\sinh w/(\cosh w-\cos t)\) is the Poisson kernel of radius \(e^{-w}\), whose \(k\)-th Fourier coefficient is \(e^{-|k|w}\). Thus the unnormalized transform is proportional to
\[
 \int_{\mathbb R}\int_{|x|}^\infty
 \frac{e^{(a-b+i\xi)x-|k|w}}{(\cosh w-\cosh x)^C}\,dw\,dx.
\]
Set \(r=(w+x)/2\) and \(s=(w-x)/2\), so the integration region is the positive quadrant and \(\cosh w-\cosh x=2\sinh r\,\sinh s\). Writing \(\mathrm B\) for the beta function, the substitutions \(e^{-2r}\) and \(e^{-2s}\) give
\(\mathrm B(b+\bar u,q)\mathrm B(a+u,q)\) times a constant independent of frequency. The integrand is bounded in modulus by the zero-frequency integrand, whose integral is finite because \(a,b,q>0\). Dividing by that zero-frequency value proves \eqref{eq:2.5}. In particular, \(\widehat p_*=|r_u(C/2)|^2>0\).

The Euler product for \(\Gamma\) gives
\[
 \log |r_u(A)|=-\int_0^\infty e^{-Av}
 \frac{(1-e^{-qv})(1-\Re e^{-uv})}{v(1-e^{-v})}\,dv .
\]
To verify the integral, expand \((1-e^{-v})^{-1}\) as a geometric series and use \(\Re\int_0^\infty(e^{-Av}-e^{-Bv})\,dv/v=\log|B/A|\) for \(A,B\) of positive real part. Subtract the expressions at \(A=a,b\) from twice the expression at \(A=C/2\). Since \(e^{-av}+e^{-bv}-2e^{-Cv/2}=(e^{-av/2}-e^{-bv/2})^2\), this gives
\[
 \log\frac{\widehat p_*(\xi,k)}{|\widehat p_{a,b}(\xi,k)|}
 =\int_0^\infty (e^{-av/2}-e^{-bv/2})^2
 \frac{(1-e^{-qv})(1-e^{-|k|v/2}\cos(\xi v/2))}{v(1-e^{-v})}\,dv .
 \tag{2.6}\label{eq:2.6}
\]
If \(a\ne b\), the integrand is positive on a set of positive measure at every nonzero frequency. Together with \(\widehat p=\widehat p_{a,b}\widehat\nu\) from \eqref{eq:2.3}, this proves the Fourier bound, strictly at nonzero frequencies when \(N>1\) or \(a\ne b\).

To conclude the proof, apply Lemma~\ref{th:2.2} to \(p\) and \(p_*\). The mass forms are continuous by the argument after \eqref{eq:2.2}, and the infimum for \(Q\) is attained. When \(Q\) is not similar to the isosceles triangle with corresponding base, the Fourier inequality is strict at \((\xi,0)\) for every \(\xi\ne0\). The lemma then gives strict inequality in \eqref{eq:2.1}. Similarity with corresponding base gives equality.\qed

\section{Convexity for small apex angles}\label{sec:3}

\begin{proposition}[Small apex angles]\label{th:3.1}
Let \(J=\lambda_1(\mathbb D)\), where \(\mathbb D\) is the unit disk. For \(0<r\le1/7\),
\[
 F''(2\pi r)>\frac{11Jr^2}{4\pi},\qquad
 F(2\pi r)\le\frac{16}{15}\pi Jr .
\]
\end{proposition}

\noindent\textit{Proof.}\enspace We map a circular sector to an isosceles triangle and estimate the area and first mixed eigenvalue of the image separately.

\subsection{A conformal map from a sector}

We first use a real parameter \(0<s<1/2\) for the conformal map. Later, \(r\) will be the center of a disk of complex parameters \(\sigma\). On the sector \(0<\arg z<2\pi s,\ |z|<1\), take the primitive of \((1-z^{1/s})^{-2s}\) that vanishes at zero. Use the branch with \(0<\arg z<2\pi s\) for \(z^{1/s}\) and the analytic power of \(1-z^{1/s}\) equal to one at zero. The bounding radii map to straight segments of equal length. On the circular arc \(z=e^{i\theta}\), \(\arg(1-e^{i\theta/s})=\theta/(2s)-\pi/2\) for \(0<\theta<2\pi s\), so the image tangent has constant argument \(\pi/2+\pi s\). The endpoint singularities are integrable. The map takes the boundary once around a triangle, so the argument principle gives a conformal bijection. This map and the weight in \eqref{eq:3.2} also occur in Wang's Appendix S1.1 \cite{ref25}.

Integration along a radius gives leg length \(s\,\mathrm B(s,1-2s)\). The area is half the square of this length times \(\sin(2\pi s)\). The Gamma reflection identity writes it as \(\pi s Z(s)\), with \(Z\) as below. Expanding the logarithm of the Gamma function, with \(\zeta(m)=\sum_{j\ge1}j^{-m}\) for \(m>1\), gives
\[
 Z(s)=\frac{\Gamma(1+s)^2\Gamma(1-2s)}{\Gamma(1-s)^2\Gamma(1+2s)},\qquad
 \log Z(s)=\sum_{\substack{m\ge3\\ m\ {\rm odd}}}
     \frac{2^{m+1}-4}{m}\zeta(m)s^m. \tag{3.1}\label{eq:3.1}
\]
Put \(\tau=-\log|z|\) and \(\varphi=\arg(z)/s\). In the following forms integrate with respect to \(d\tau\,d\varphi/(2\pi)\) on \((0,\infty)\times(\mathbb R/2\pi\mathbb Z)\), and take the forms linear in the first argument. The two radii have the same trace for the mixed ground state, by reflection in the bisector, so we use periodic functions of \(\varphi\). Conversely, cutting a periodic test along the identified radii gives an admissible function for the mixed problem with the same quotient. The apex has zero capacity. For the energy \(E_s\) and mass \(m_s\) below, write \(E_s(u)=E_s(u,u)\). Let \(L(s)\) be the minimum of \(E_s(u)/m_s(u,u)\) over nonzero periodic functions of finite energy with zero trace at \(\tau=0\). Conformal invariance and the area formula give
\[
 \begin{aligned}
 E_s(u,v)&=\int(u_\tau\overline{v_\tau}+s^{-2}u_\varphi\overline{v_\varphi}),\\
 m_s(u,v)&=\int e^{-2\tau}|1-e^{-\tau/s+i\varphi}|^{-4s}u\bar v,\\
 F(2\pi s)&=\pi sZ(s)L(s),
 \end{aligned} \tag{3.2}\label{eq:3.2}
\]
For each fixed \(r>0\), use the norm \(E_r^{1/2}\) on this space, including the zero function. The zero trace and Cauchy–Schwarz give \(\|u(\tau,\cdot)\|_2^2\le\tau E_r(u)\), with angular measure \(d\varphi/(2\pi)\). This estimate controls convergence on every finite cylinder and proves completeness. Write
\[
 \frac{F(2\pi r)}{\pi}=JrZ(r)+rZ(r)(L(r)-J).
\]
We bound the second derivative of the first term below by \(54Jr^2\) and the absolute value of the second derivative of the remainder above by \(43Jr^2\).

\subsection{The disk ground state and the radial gap}

Set \(m_0(u,v)=\int e^{-2\tau}u\bar v\). Rellich compactness on bounded cylinders and the slice estimate at infinity give compactness for this mass. Under the substitution \(\rho=e^{-\tau}\), the quotient \(E_r(u)/m_0(u,u)\) for functions independent of \(\varphi\) is the radial Dirichlet quotient on the unit disk. The disk ground state is radial by rotational invariance and simplicity. Denote it in the \(\tau\)-coordinate by \(B\) and normalize it by \(m_0(B,B)=1\). The trial function \(1-e^{-2\tau}\) has energy one and mass \(1/6\), so \(0<J\le6\). The radial equation is \(B''=-Je^{-2\tau}B\). Since \(B\) is positive and concave and \(B'\) is square integrable, \(B'(\tau)\to0\) as \(\tau\to\infty\). Multiply the equation by \(2B'\) and integrate. The boundary terms \(e^{-2\tau}B^2\) vanish at both ends by the trace and the bound \(|B(\tau)|^2\le\tau E_r(B)\). This gives
\[
 B'(0)^2=2J,\qquad 0\le B(\tau)\le\sqrt{2J}\,\tau. \tag{3.3}\label{eq:3.3}
\]
Averaging any test function over the circle in the radial equation gives \(E_r(B,v)=Jm_0(B,v)\).

Let \(V\) be the subspace of functions with \(m_0(v,B)=0\). By the min–max principle and zero extension, the second radial eigenvalue is at least the second Dirichlet eigenvalue of the disk, which is at least that of \((-1,1)^2\). The latter is \(5\pi^2/4\) by separation of variables. Thus for radial \(v\in V\),
\[
                   m_0(v,v)\le\frac{4}{5\pi^2}E_r(v).
\]

\subsection{Fourier estimates for the mass weight}

Fix \(0<r\le1/7\). For a complex parameter \(|\sigma-r|\le r/3\), put \(c=3/4\), \(x=c\tau/r\) and \(p=|\sigma|\). Then \(p\le r/c\). The quotient \(r/\sigma\) lies in the disk with real endpoints \(3/4\) and \(3/2\). On this disk the real part is at least \(c\) and the real part of the square is at least \(c^2\). Indeed, for real part \(y\), the squared imaginary part is at most \((y-3/4)(3/2-y)\), which is at most \(y^2-(3/4)^2\). Define \(E_\sigma\) by using \(\sigma^{-2}\) in the angular term of \eqref{eq:3.2}. Define \(m_\sigma\) by replacing the factor after \(e^{-2\tau}\) by
\[
 W_\sigma=(1-e^{-\tau/\sigma+i\varphi})^{-2\sigma}
          (1-e^{-\tau/\sigma-i\varphi})^{-2\sigma}.
\]
Use the power series for \(\log(1-w)\), \(|w|<1\), to define the factors. Since \(\Re(1/\sigma)\ge c/r\), the binomial series gives \(|W_\sigma|\le(1-e^{-x})^{-4p}\), with \(4p\le16/21\). For test functions of bounded energy, the slice estimate then bounds the mass on \(0<\tau<\epsilon\) by \(O_r(\epsilon^{26/21})\) and that on \(\tau>R\) by \(O_r((1+R)e^{-2R})\). On each truncated cylinder the series converge locally uniformly. Together these bounds give holomorphic dependence of \(m_\sigma\) in operator norm. For real parameters, Rellich compactness on the truncated cylinders and the same bounds at the ends give compactness of the mass embedding.

Let \(\langle\cdot\rangle\) denote averaging in \(\varphi\) with measure \(d\varphi/(2\pi)\), and put \(\mathcal V=\langle|W_\sigma-\langle W_\sigma\rangle|^2\rangle\). The binomial expansion gives
\[
 \begin{aligned}
 |W_\sigma-1|&\le4p\sum_{n\ge1}e^{-nx},\\
 |\langle W_\sigma\rangle-1|&\le4p^2\sum_{n\ge1}e^{-2nx},\\
 \mathcal V&\le16p^2\left[\tfrac12 e^{-2x}+\sum_{n\ge2}e^{-2nx}\right].
 \end{aligned} \tag{3.4}\label{eq:3.4}
\]
The binomial coefficients are \(a_0(v)=1\) and \(a_n(v)=v(v+1)\cdots(v+n-1)/n!\) for \(n\ge1\). They satisfy \(|a_n(v)|\le a_n(|v|)\le|v|\) if \(n\ge1\) and \(|v|\le1\). To bound the modulus of the product, replace both factors by \((1-e^{-x})^{-2p}\). Their product has coefficients \(a_n(4p)\), giving the first bound. For the mean only equal powers from the two factors contribute, giving the coefficients \(a_n(2p)^2\) in the second bound. For the variance put \(\rho=e^{-x}\) and \(H_\rho=|1-\rho e^{i\varphi}|^{-4p}\). Each Fourier coefficient of \(W_\sigma\) is bounded in modulus by the corresponding nonnegative coefficient of \(H_\rho\). Parseval's identity gives
\[
 \mathcal V\le\langle H_\rho^2\rangle-\langle H_\rho\rangle^2,\qquad
 \langle H_\rho^2\rangle=\sum_{j\ge0}a_j(4p)^2\rho^{2j},\qquad
 \langle H_\rho\rangle=\sum_{j\ge0}a_j(2p)^2\rho^{2j}.
\]
The coefficient of \(\rho^2\) in the difference is \((4p)^2-2(2p)^2=8p^2\). For every power \(\rho^{2j}\) with \(j\ge2\), the coefficient is at most \((4p)^2\). This proves the last line of \eqref{eq:3.4}.

\subsection{Bounds for the mass perturbation}

Set \(\Delta=m_\sigma-m_0\), and take operator and dual norms with respect to \(E_r^{1/2}\). For \(j>1\) and \(t\ge0\), set \(S_j(t)=\sum_{n\ge1}(n+t)^{-j}\). By convexity, \((n+t)^{-j}\) is bounded by the integral of \((x+t)^{-j}\) over \([n-\tfrac12,n+\tfrac12]\). Summing this bound for \(n\ge2\) gives
\[
 S_j(t)\le U_j(t):=(1+t)^{-j}+\frac{(3/2+t)^{1-j}}{j-1}.
\]
For fixed \(r\), abbreviate \(D(r)\) and \(H(r)\) to \(D\) and \(H\), where
\[
 D(r)=\frac{2}{c^5}S_3(r/c),\qquad
 H(r)=\frac{8}{c^5}\left[S_3(2r/c)-\tfrac12(1+2r/c)^{-3}\right].
\]
Then
\[
 \begin{aligned}
 \|\Delta\|&\le 4r^3c^{-3}S_2(2r/c),\qquad
 |\Delta(B,B)|\le Jr^5D,\\
 \max\{\|\Delta(B,\cdot)\|^2,\|\Delta(\cdot,B)\|^2\}
 &\le Jr^7[H+(16/3)r^2D].
 \end{aligned} \tag{3.5}\label{eq:3.5}
\]
For the operator bound, apply the first line of \eqref{eq:3.4} and the slice inequality in both arguments, then integrate \(\tau e^{-(2+nc/r)\tau}\). For two radial functions, the second line of \eqref{eq:3.4} bounds the absolute value of the difference of weights by the density of the measure
\[
 d\omega=4p^2 e^{-2\tau}\sum_{n\ge1}e^{-2nc\tau/r}\,d\tau .
\]
By \eqref{eq:3.3}, \(\int B^2d\omega\le Jr^5D\). For a radial test function \(v\), the slice bound gives
\[
 \int|v|^2d\omega\le r^4c^{-4}S_2(r/c)E_r(v)\le\frac{16}{3}r^4E_r(v),
\]
since \(c^{-4}U_2(0)=1280/243<16/3\). Cauchy–Schwarz for \(d\omega\) bounds the squared dual norm on radial test functions by \((16/3)Jr^9D\). For a test function of angular mean zero, the angular Poincaré inequality gives \(\int |v|^2\le r^2E_r(v)\). Cauchy–Schwarz on the circle, after subtracting the mean of \(W_\sigma\), bounds the squared dual norm by \(2Jr^2\int_0^\infty\tau^2e^{-4\tau}\mathcal V\,d\tau\le Jr^7H\). The squared bounds add because the two subspaces are orthogonal for \(E_r\). These integrations use \(\int_0^\infty \tau^j e^{-d\tau}\,d\tau=j!/d^{j+1}\) for integers \(j\ge0\) and \(d>0\).

\subsection{A second-derivative estimate}

For \(|\lambda-J|\le1/4\), write \(\Lambda=25/4\), so \(|\lambda|\le\Lambda\). On \(V\), the radial gap and the angular Poincaré inequality give a lower bound \(\min(1-4\Lambda/(5\pi^2),c^2-\Lambda r^2)E_r\) for the real part of \(E_\sigma-\lambda m_0\). The bound on \(\|\Delta\|\) in \eqref{eq:3.5} increases with \(r\) term by term. At \(r=1/7\), \(U_2(8/21)<17/16\). Using this upper bound on \(|\lambda|\|\Delta\|\) gives, for \(v\in V\),
\[
 \Re(E_\sigma-\lambda m_\sigma)(v,v)\ge
 \left[\frac{341}{784}-\frac{1700}{9261}\right]E_r(v)
 >\frac14 E_r(v)\qquad(v\ne0). \tag{3.6}\label{eq:3.6}
\]
Let \(\mathcal A:V\to V^*\) be the operator induced on \(V\) by \(E_\sigma-\lambda m_\sigma\), where \(V^*\) is its continuous conjugate dual. By \eqref{eq:3.6} and the complex Lax–Milgram theorem, \(\|\mathcal A^{-1}\|\le4\). The inverse depends holomorphically on \(\sigma,\lambda\) by a local Neumann series.

By \eqref{eq:3.6}, an eigenfunction with \(|\lambda-J|\le1/4\) has a nonzero component along \(B\). Normalize it to \(B+v\), with \(v\in V\). Testing on \(V\) gives \(v=\lambda\mathcal A^{-1}\Delta(B,\cdot)\). Testing against \(B\) gives
\[
 \lambda[1+\Delta(B,B)+\Delta(v,B)]-J=0 .
\]
By \eqref{eq:3.5} its left side differs from \(\lambda-J\) by at most \(Jr^5Q(r)\), where
\[
 Q(r)=\Lambda D+4\Lambda^2\{r^2H+(16/3)r^4D\}.
\]
The product \(r^2Q(r)\) is increasing for \(0<r\le1/7\). Indeed, \(r^2/(n+r/c)^3\) increases on this interval for \(n\ge1\), and \(r^4/(n+2r/c)^3\) increases for all \(r>0\). All the series coefficients in \(H\) are positive. The bound by \(U_3\) gives
\[
 D(1/7)\le\frac{2}{c^5}U_3(4/21)<\frac{13}{2},\qquad
 H(1/7)\le\frac{8}{c^5}\left[U_3(8/21)-\tfrac12(21/29)^3\right]<12 .
\]
Hence \(49r^2Q(r)\le Q(1/7)<\Lambda(13/2)+4\Lambda^2(12/49+35/49^2)<82\), and \(Jr^5Q(r)<6\cdot82/7^5<1/4\). Rouché's theorem on \(|\lambda-J|=1/4\) gives one root, counted with multiplicity. Denote it by \(\mathcal L(\sigma)\). It is simple and hence holomorphic in \(\sigma\), with \(|\mathcal L(\sigma)-J|\le Jr^5Q(r)\).

For real \(\sigma\) in the disk, uniqueness makes the root real and the equations on \(V\) and \(B\) give an eigenfunction. If the second eigenvalue were at most \(J+1/4\), the span of the first two eigenfunctions would meet \(V\) in a nonzero vector with \((E_\sigma-(J+1/4)m_\sigma)(v,v)\le0\). This contradicts \eqref{eq:3.6}. Thus \(\mathcal L(\sigma)=L(\sigma)\) for real \(\sigma\) in the disk.

The positive coefficients in \eqref{eq:3.1}, \(\zeta(3)<5/4\) and \(|\sigma|\le4/21\) give
\[
 |\log Z(\sigma)|\le\frac{20}{3}\frac{|\sigma|^3}{1-4|\sigma|^2}<\frac{1}{16},\qquad |Z(\sigma)|<\frac{16}{15} .
\]
On \(|\sigma-r|\le r/3\), the holomorphic remainder \(\sigma Z(\sigma)(\mathcal L(\sigma)-J)\) has modulus at most \((4r/3)(16/15)Jr^5Q(r)\). Cauchy's estimate gives
\[
 \bigl|[rZ(r)(L(r)-J)]''\bigr|
 \le\frac{2}{(r/3)^2}\frac{4r}{3}\frac{16}{15}Jr^5Q(r)<43Jr^2 .
\]
The cubic coefficient of \(Z\) is \(4\zeta(3)\), and all its coefficients are nonnegative. Hence \(J[rZ(r)]''\ge48\zeta(3)Jr^2>54Jr^2\), since \(\zeta(3)>9/8\). Combining the two second-derivative estimates gives \([rZ(r)L(r)]''>11Jr^2\). By \eqref{eq:3.2} and the chain rule this gives the asserted bound on \(F''(2\pi r)\). For the upper bound, the circular mean of \(|1-e^{-\tau/r+i\varphi}|^{-4r}\) equals \(\sum_{n\ge0}a_n(2r)^2e^{-2n\tau/r}\), which is at least its constant term one. Testing with the radial state \(B\) gives \(L(r)\le J\). Now use \(Z(r)<16/15\) in \eqref{eq:3.2}.\qed

\section{Supporting lines for the triangle eigenvalue}\label{sec:4}

\begin{proposition}[Supporting-line inequality]\label{th:4.1}
The function \(F\) is strictly convex on \((0,\pi/2]\). For each \(0<a\le\pi/2\), its tangent at \(a\) lies strictly below \(F(\alpha)\) at every \(\alpha\in[0,\pi)\) other than \(a\). In particular \(aF'(a)>F(a)>0\).
\end{proposition}

\subsection{A fixed triangle and a trial function}

\noindent\textit{Proof.}\enspace We first prove strict convexity for apex angles at most \(\pi/2\). We then compare \(F\) with its tangent at \(\pi/2\) for obtuse angles; this extends the supporting-line inequality to the whole interval.

For \(0<\alpha<\pi\), set \(t=\tan(\alpha/2)\) and \(z=t^{-2}\). The triangle \(I_\alpha\) has vertices \((0,0),(1,\pm t)\), area \(t\), and an even ground state across the horizontal axis. The map \((x,y)\mapsto(x,ty)\) sends \(D=\{(x,y):0<y<x<1\}\) onto its upper half. On \(H^1(D)\) with zero trace on \(x=1\) define
\[
 A(u,v)=\int_Du_x\overline{v_x}\,d\mu,\qquad
 B(u,v)=\int_Du_y\overline{v_y}\,d\mu,\qquad
 m(u,v)=\int_Du\bar v\,d\mu,\qquad d\mu=2\,dx\,dy .
\]
Write \(A(v)=A(v,v)\), \(B(v)=B(v,v)\), and \(\|v\|_2^2=m(v,v)\). Orthogonality in this section is for \(m\). Let \(G(z)\) be the first eigenvalue for the energy \(Q_z=A+zB\) and mass \(m\). Changing variables in the Rayleigh quotient gives
\[
 F(\alpha)=tG(z).
\]
The Rayleigh principle makes \(G\) increasing and concave. For \(z>0\), the ground state is simple and its eigenvalue is analytic. Indeed, the inverse after a fixed positive mass shift depends locally analytically on \(z\) by a Neumann series in the energy norm. Its compactness and a contour around the isolated eigenvalue give an analytic projection onto the one-dimensional eigenspace.

At \(z=1\), even reflection across the diagonal extends the problem to the symmetric subspace on the unit square, with Neumann conditions on the axes and Dirichlet conditions on \(x=1\) and \(y=1\). The normalized ground state on \(D\) is \(e=2\cos(\pi x/2)\cos(\pi y/2)\), with \(\lambda=G(1)=\pi^2/2\). Differentiation gives \(G'(1)=B(e,e)\), so integration gives \(\eta=G'(1)=\lambda/2-2\). The next sum of two odd squares after \(1+1\) is \(1+9\), so separation of variables gives \(Q_1(w,w)\ge5\lambda\|w\|_2^2\) for \(w\perp e\). For fixed \(y\), the Dirichlet–Neumann inequality on \(y<x<1\) gives \(A(w)\ge(\lambda/2)\|w\|_2^2\) after integrating in \(y\).

For \(1\le z\le5\), use the trial function \(v=(1-x^2)(3-x^2-y^2)\). The identity \(\int_D x^i y^j\,dx\,dy=1/[(1+j)(2+i+j)]\) for integers \(i,j\ge0\) gives \(A(v)=569/45\), \(B(v)=1/9\), and \(m(v,v)=539/225\). Put \(\delta=z-1\), \(q=G(z)-\lambda\), and \(h=5\lambda-G(z)\). The trial quotient gives
\[
 p(\delta)=\frac{5(574+5\delta)}{539}-\lambda,\qquad
 0\le q\le p(\delta),\qquad h>19 . \tag{4.1}\label{eq:4.1}
\]
Here and below we use \(223/71<\pi<22/7\), so in particular \(\lambda>4.93\). All terminating decimals denote exact rationals. The bounds on \(\pi\) follow by using four and five terms of the alternating series for \(\arctan(1/5)\), and one and two for \(\arctan(1/239)\), in \(\pi/4=4\arctan(1/5)-\arctan(1/239)\). The last identity follows from the tangent addition formula.

\subsection{Convexity for the remaining acute angles}

We project the ground state at \(z\) onto the orthogonal complement of the ground state at \(z=1\), as in Wang's (S1.7)–(S1.8) \cite{ref25}. Throughout this calculation \(1\le z\le5\). We abbreviate \(G(z)\) to \(G\) and write \(R=z^2(z+1)\) and \(g=G'(z)\). Concavity gives \(0\le g\le\eta<0.47\). On \(e^\perp\) use the inner product \(\langle v,w\rangle_H=Q_1(v,w)-G(z)m(v,w)\), which is positive definite by \eqref{eq:4.1} and the gap at \(z=1\). For each fixed \(z\), define the operators \(T,\mathcal M\) and vector \(f\) there by \(\langle Tv,w\rangle_H=B(v,w)\), \(\langle\mathcal Mv,w\rangle_H=m(v,w)\), and \(\langle f,w\rangle_H=B(e,w)\). Take norms with respect to this inner product, and let \(I\) denote the identity operator. Then
\[
 0\le T\le\kappa I,\quad \kappa=\frac{9\lambda/2}{5\lambda-G}<1.2,\qquad
 \|\mathcal M\|\le1/h,\qquad f=T^{1/2}\xi,\quad\|\xi\|_H^2\le\eta .
\]
Indeed, \(B\le Q_1-(\lambda/2)m\) and \((v-\lambda/2)/(v-G(z))\) decreases for \(v\ge5\lambda\). The gap also bounds \(\|\mathcal M\|\) by \(1/h\). Cauchy–Schwarz for \(B\) gives \(|\langle f,w\rangle_H|\le\sqrt{\eta}\|T^{1/2}w\|_H\). The Riesz representation theorem on the closure of the range of \(T^{1/2}\) gives the stated vector \(\xi\).

On \(e^\perp\), let \(C=(I+\delta T)^{-1}\), and write the normalized positive ground state as \(u=c(e-y)\) with \(y\perp e\) and \(c=m(u,e)>0\). Testing the eigenvalue equation on \(e^\perp\) and then on \(e\) gives
\[
 u=c(e-y),\qquad y=\delta Cf,\qquad c=(1+\|y\|_2^2)^{-1/2},
 \qquad q=\delta\eta-\delta^2\langle f,Cf\rangle_H. \tag{4.2}\label{eq:4.2}
\]
Set \(\mathcal S_z=Q_z-G(z)m\). We differentiate the variational equation before expressing it through these \(z\)-dependent operators. For fixed \(w\perp e\), differentiate \(\mathcal S_z(y,w)=\delta B(e,w)\) in \(z\). This gives
\[
 \mathcal S_z(y',w)=B(e-y,w)+g\,m(y,w).
\]
Now use the inner product at the current \(z\), together with \(y=\delta Cf\) and \(I-\delta TC=C\), to obtain
\[
 C^{-1}y'=f-Ty+g\mathcal M y=Cf+\delta g\mathcal M Cf,\qquad
 y'=C^2f+\delta g C\mathcal M Cf .
\]
For the real normalized ground state, differentiation of \(m(u,u)=1\) gives \(m(u,u')=0\), and differentiation of the eigenvalue equation gives \(G'=B(u,u)\). Testing the differentiated equation on \(u'\) gives \(-G''/2=\mathcal S_z(u',u')\). Write \(K=-G''/2\). Since \(u'=(c'/c)u-cy'\) and \(\mathcal S_z(u,\cdot)=0\), we have \(K=c^2\|C^{-1/2}y'\|_H^2\). Put \(K_0=\langle f,C^3f\rangle_H\). For \(0\le v\le6/5\) and \(0\le\delta\le4\),
\[
 3(1+\delta v)^2-z(z+1)v
 =3-2v+(3-\delta)\delta v+3(\delta v)^2
 \ge3/5-\delta v+3(\delta v)^2>0 .
\]
Thus \(Rv/(1+\delta v)^2\le3z\) on the spectrum of \(T\). Since \(f=T^{1/2}\xi\), the spectral theorem and \(0<C\le I\) give \(RK_0\le R\|Cf\|_H^2\le3z\eta\). The formula for \(y'\), together with \(c\le1\), \(\|C^{1/2}\|\le1\) and \(\|\mathcal M\|\le1/h\), now gives
\[
 \sqrt{RK}\le\sqrt{RK_0}+\frac{\delta g}{h}\sqrt{3z\eta}.
 \tag{4.3}\label{eq:4.3}
\]
Since \(dt/d\alpha=(1+t^2)/2\) and \(dz/d\alpha=-z(1+t^2)/t\), we have \(F'=(1+t^2)(G-2zg)/2\). With \(M=2F''(\alpha)/(t(1+t^2))\), the second derivative gives
\[
 M=G-2zg+z(z+1)(g+2zG'')=G+z\delta g-4RK .
\]
Squaring \eqref{eq:4.3} and using \(RK_0\le3z\eta\) and \(g\le\eta\), we find that \(4RK-4RK_0\) is at most \(z\delta g(24\eta/h+12\delta\eta^2/h^2)\). The factor in parentheses is less than one because \(\eta<0.47\), \(h>19\), and \(\delta\le4\). Thus
\[
                   M\ge G-4RK_0\qquad(1\le z\le5). \tag{4.4}\label{eq:4.4}
\]

The following two bounds give \(F''>0\) up to the range covered by Proposition~\ref{th:3.1}. They will also be used in the proof of Theorem~\ref{th:6.2}.
\[
 M>5/7\quad(1\le z\le5/2),\qquad
 M>21/50\quad(5/2\le z\le9/2). \tag{4.5}\label{eq:4.5}
\]
For the first range, the spectrum of \(T\) lies in \([0,\kappa]\), where \(\kappa=9\lambda/(2h)\). For \(1\le v\le\kappa\),
\[
 \frac{Rv}{(1+\delta v)^3}\le\left(1+\frac1z\right)\kappa\le\frac{9}{4}.
\]
Indeed, concavity gives \(h\ge4\lambda-\delta\eta\ge2\lambda(1+1/z)\). For \(z>1\) the second inequality follows from \(z\eta\le2\lambda\), using \(z\le5/2\), \(\eta<0.47\), and \(\lambda>4.93\). At \(z=1\) there is equality in the second inequality. On \(0\le v\le1\) and \(1\le z\le3/2\), the maximum is at \(v=1\) and equals \(1+1/z\le2\). For \(3/2\le z\le5/2\), the maximum over \(v\ge0\) is \(4R/(27\delta)\). The function
\[
 R/\delta=\delta^2+4\delta+5+2/\delta
\]
is convex on \([1/2,3/2]\), so its endpoint values give \(4R/(27\delta)\le175/81<9/4\). The spectral theorem now gives \(4RK_0\le9\eta\), and \eqref{eq:4.4} yields
\[
 M\ge\lambda-9\eta=18-7\pi^2/4>5/7.
\]

For the second range, \(3/2\le\delta\le7/2\). The operator inequality \((I-\delta T)^2\ge0\) gives \(4\delta TC^3\le C\). Equation \eqref{eq:4.2} and \(\|\xi\|_H^2\le\eta\) imply
\[
 4\delta^2K_0\le\delta\langle\xi,C\xi\rangle_H\le q, \tag{4.6}\label{eq:4.6}
\]
since \(q/\delta=\eta-\|\xi\|_H^2+\langle\xi,C\xi\rangle_H\). Consequently
\[
 M\ge\lambda-(R/\delta^2-1)q
 \ge\lambda-(\delta+3+5/\delta+2/\delta^2)p(\delta).
\]
The rightmost expression is concave in \(\delta\). Indeed \(p(\delta)\) is affine with positive constant and linear coefficients, so the product being subtracted is a positive linear combination of \(\delta^{-2},\delta^{-1},1,\delta,\delta^2\). Replacing \(\lambda\) by \(493/100\) decreases the expression because the coefficient of \(\lambda\) is \(1+\delta+3+5/\delta+2/\delta^2>0\). At both \(\delta=3/2\) and \(\delta=7/2\), substitution gives a value greater than \(21/50\). This proves the second bound in \eqref{eq:4.5}.

To verify the overlap with Proposition~\ref{th:3.1}, use \(\pi/7>4/9\) and integrate \((\tan u)'=1+\tan^2u>1+u^2\) for \(0<u<\pi/2\). This gives
\[
 \tan(\pi/7)>\frac49+\frac13\left(\frac49\right)^3
 =\frac{1036}{2187}>\frac{\sqrt2}{3}.
\]
Hence \(\cot^2(\pi/7)<9/2\), completing the curvature proof on \((0,\pi/2]\).

\subsection{The tangent at a right apex angle}

Since \(G(1)=\lambda\) and \(G'(1)=\lambda/2-2\), we have \(F(\pi/2)=\lambda\) and \(F'(\pi/2)=4\). For \(t>1\) write \(d=(t-1)/(t+1)\). Then \(\alpha-\pi/2=2\arctan d\le2d\). To show that the tangent lies strictly below \(F(\alpha)\), it is enough to prove \(F(\alpha)>J_t\), where \(J_t=\lambda+8d\).

For \(w\perp e\), integrate by parts in \(y\). The interior term is \((\lambda/2)m(e,w)=0\). Using \(w(1,y)=0\) to write the remaining boundary value as an integral of \(w_x\) gives
\[
 B(e,w)=-\pi\int_0^1\sin(\pi y)\overline{w(y,y)}\,dy
        =\int_D \tfrac{\pi}{2}\sin(\pi y)\overline{w_x}\,d\mu .
\]
The identity \(2\int_0^1(1-y)\sin^2(\pi y)\,dy=1/2\) and Cauchy–Schwarz give the first estimate below. Applying Cauchy–Schwarz to \(B\) gives the second.
\[
 |B(e,w)|^2\le(\lambda/4)A(w),\qquad |B(e,w)|^2\le\eta B(w).
\]

For \(1<t\le2\), put \(\rho=J_t/(5\lambda)=1/5+8d/(5\lambda)<1/3\), and let \(S=4(t-\rho)/\lambda+(t^{-1}-\rho)/\eta>0\). The gap on \(e^\perp\) and the two bounds on \(B(e,w)\) give
\[
 (tA+t^{-1}B-J_tm)(w,w)\ge (t-\rho)A(w)+(t^{-1}-\rho)B(w)\ge S|B(e,w)|^2,
\]
Decompose the normalized positive ground state at \(t\) as \(ce+w\) with \(c>0\) and \(w\perp e\). Since \(A(e,w)=-B(e,w)\), the cross term is \(2c(t^{-1}-t)\Re B(e,w)\). Completing the square gives
\[
 F(2\arctan t)-J_t\ge c^2\left[tA(e)+B(e)/t-J_t-\frac{(t^{-1}-t)^2}{S}\right]. \tag{4.7}\label{eq:4.7}
\]
Put \(W(t)=t\lambda S/2=2t^2-2\rho t+\lambda(1-\rho t)/(2\eta)\).

The term \(tA(e)+B(e)/t-J_t\) equals \((t-1)^2(\lambda/2+2d)/t\). Cancelling \((t-1)^2\) in the bracket of \eqref{eq:4.7}, we need \((1+4d/\lambda)W(t)>(t+1)^2\). As \(\rho t<1\), \(\rho<1/3\), and \(\lambda/(2\eta)>5\), we have \(W(t)\ge 2t^2-7t/3+5\ge14/3\). Since \(t+1\le3\), this gives \(4W(t)/(t+1)\ge56/9>\lambda\). Hence
\[
 \begin{aligned}
 (1+4d/\lambda)W(t)-(t+1)^2
 &>W(t)+(t-1)-(t+1)^2\\
 &\ge 2t^2-7t/3+5+(t-1)-(t+1)^2\\
 &=(t-5/3)^2+2/9>0.
 \end{aligned}
\]

For \(t\ge2\), use concavity of \(G\) between \(z=0\) and \(z=1/3\). Since \(G\) is increasing and the horizontal Dirichlet–Neumann inequality gives \(G(z)\ge\lambda/2\), the limit \(G(0)=\lim_{z\downarrow0}G(z)\) exists and is at least \(\lambda/2\). At \(z=1/3\), the mixed triangle joins a side of an equilateral triangle of inradius one to its center. Let \(d_1,d_2,d_3\) be the inward distances to the three sides of the equilateral triangle. They sum to three. The function \(\sum_i\sin(2\pi d_i/3)=4\prod_i\sin(\pi d_i/3)\) is positive inside, vanishes on the boundary, and solves the eigenvalue equation with eigenvalue \(4\pi^2/9\). By reflection, \(G(1/3)=4\pi^2/9\). The chord bound gives
\[
                    F(\alpha)\ge\frac{\pi^2}{12}\left(3t+\frac{7}{t}\right).
\]
The comparison with \(J_t\) follows from the exact identity
\[
 \frac{\pi^2}{12}(3t+7/t)-\frac{\pi^2}{2}-\frac{\pi^2}{7}t
       =\frac{\pi^2(3t-7)^2}{84t}.
\]
Indeed, \(t/d=3+(t-1)+2/(t-1)\ge3+2\sqrt2\), so \(8d/t\le8/(3+2\sqrt2)<7/5<\pi^2/7\). Thus \(F(\alpha)>J_t\) also for \(t\ge2\).

Strict convexity gives the tangent inequality on \((0,\pi/2]\), with equality only at the point of tangency. For \(a<\pi/2\), the tangent at \(a\) has smaller slope than the tangent at \(\pi/2\), and has smaller value there. Thus it also lies below the graph for \(\alpha>\pi/2\). Since \(F(0+)=0\) and \(F>0\), strict convexity gives \(F'(a)>F(a)/a>0\). Hence the tangent has negative value at zero.\qed

\section{A partition of a simple polygon}\label{sec:5}

The partition must make the supporting line give the same mixed eigenvalue bound on every cell. If that line is \(s(\gamma-b)\), with \(s,b>0\), then \(\gamma\ge b+c|V|\) gives \(\lambda_1^{\mathrm{DN}}(V)\ge sc\) for a cell \(V\), by \eqref{eq:1.2}. Here \(c>0\), and \(\gamma\) is \(\pi\) minus the sum of the base angles. Splitting a cell preserves its total area but adds one copy of \(b\) to the sum of the required angle bounds. We therefore assign an additional angle \(b\) to each vertex through which a cut is needed. The construction has three steps: choose weighted cells in the convex hull, cut through the remaining polygon vertices, and select the pieces inside the polygon. The last step may require changing a cell's base to a side on the polygon boundary.

\Needspace{8\baselineskip}
\begin{lemma}[Partition into convex cells]\label{th:5.1}
Let \(P\) be a simple polygon of positive area with \(k\ge3\) sides. Suppose
\[
 b>0,\qquad (k-1)b<\pi,\qquad c=\frac{2\pi-kb}{|P|}. \tag{5.1}\label{eq:5.1}
\]
Then \(c>0\), and \(P\) has a finite partition, up to the cell boundaries, into positive-area convex polygons \(V_j\). Each has a whole side on \(\partial P\). With this side as the base, let \(\gamma_j\) be \(\pi\) minus the sum of the two adjacent interior angles. Then
\[
 \gamma_j\ge b+c|V_j|>0. \tag{5.2}\label{eq:5.2}
\]
\end{lemma}

\noindent\textit{Proof.}\enspace Since \(kb<k\pi/(k-1)<2\pi\), we have \(c>0\). Write the sides of the convex hull of \(\overline P\) as \(E_1,\ldots,E_m\) in cyclic order. Let \(\beta_i\in(0,\pi)\) be the hull angle at the initial endpoint of \(E_i\) and \(h_i(x)\) the distance from \(x\) to the line through \(E_i\). For \(t\in\mathbb R^m\), let \(Q_i(t)\) consist of the points in the interior of the hull at which \(e^{t_i}h_i\) is strictly smaller than every other \(e^{t_\ell}h_\ell\). These are convex cells since the distances are affine in the hull. The closure of each cell contains \(E_i\) as its base, and the cells cover the interior of the hull apart from finitely many line segments.

All indices in this construction are cyclic modulo \(m\). Let \(x_i\) be the angle between \(E_i\) and the ray from its initial endpoint on which the weighted distances to \(E_{i-1}\) and \(E_i\) agree. Equality of the two distances gives \(e^{t_i-t_{i-1}}\sin x_i=\sin(\beta_i-x_i)\). Thus
\[
 \varphi_i(v)=\arg(e^v+e^{\mathrm i\beta_i})\in(0,\beta_i),\qquad
 x_i=\varphi_i(t_i-t_{i-1}) .
\]
The corresponding two angles with the base \(E_i\) are \(x_i\) and \(\beta_{i+1}-x_{i+1}\). Their sum subtracted from \(\pi\) gives
\[
 \alpha_i=\pi-\beta_{i+1}-x_i+x_{i+1}<\pi,
 \qquad \sum_i\alpha_i=2\pi.
\]
If \(\alpha_i>0\), the two rays and \(E_i\) enclose \(Q_i(t)\) in a triangle with apex angle \(\alpha_i\). Set
\[
 \mathcal A(t)=\sum_i(\pi-\beta_{i+1})t_i
       -\sum_i\int_0^{t_i-t_{i-1}}\varphi_i(v)\,dv.
\]
Since \(\varphi_i\) decreases, \(\mathcal A\) is convex and \(\partial_i\mathcal A=\alpha_i\).

For each vertex \(X\) of \(P\) that is not a hull vertex, define \(q_X(t)\) by
\[
 q_X(t)=
 \begin{cases}
  \min_i(t_i+\log h_i(X)),&X\text{ strictly inside the hull},\\
  t_i,&X\text{ in the relative interior of }E_i .
 \end{cases}
\]
Measure distances in a fixed unit in these logarithms. We minimize the continuous convex function
\[
 \mathcal A(t)-c\int_P\min_i(t_i+\log h_i(x))\,dx
              -b\sum_i t_i-b\sum_Xq_X(t). \tag{5.3}\label{eq:5.3}
\]
The sum over \(X\) runs over the \(k-m\) vertices of \(P\) that are not hull vertices. The logarithms are integrable with respect to area. Adding a common constant to the \(t_i\) leaves \eqref{eq:5.3} unchanged because \(c|P|+kb=2\pi\).

Normalize \(\min_i t_i=0\) and put \(T=\max_i t_i\). The function \(\varphi_i\) approaches zero exponentially at \(+\infty\) and \(\beta_i\) exponentially at \(-\infty\). Integrating gives
\[
 \mathcal A(t)=\sum_i(\pi-\beta_{i+1})t_i
       +\sum_i\beta_i(t_{i-1}-t_i)_++O(1),
\]
where \(v_+=\max(v,0)\) and the error is bounded independently of \(t\). To estimate the two sums, write \(t_i\) as the integral of the indicator of \(0<h<t_i\) in \(h\), and \((t_{i-1}-t_i)_+\) as the integral of the indicator of \(t_i<h<t_{i-1}\). At each level \(0<h<T\), group the indices with \(t_i>h\) into maximal consecutive blocks in cyclic order. A block \(L,\ldots,R\) contributes
\[
 \sum_{i=L}^R(\pi-\beta_{i+1})+\beta_{R+1}
       =\pi+\sum_{i=L+1}^R(\pi-\beta_i)\ge\pi.
\]
Integrating over \(h\) gives \(\mathcal A(t)\ge\pi T-O(1)\). Using an index with \(t_i=0\) in each minimum gives uniform upper bounds on the area integral and on the \(q_X\) for vertices strictly inside the hull. The sum of the remaining terms multiplied by \(b\) is at most \(b(k-1)T\), since each \(t_i\le T\) and at least one \(t_i=0\). By \eqref{eq:5.1}, the expression to be minimized tends to infinity as \(T\to\infty\). It therefore attains a minimum on \(\min_i t_i=0\) and hence on \(\mathbb R^m\).

At a minimizing \(t\), write \(Q_i=Q_i(t)\). The partial derivative of the integral in \eqref{eq:5.3} with respect to \(t_i\) is \(|P\cap Q_i|\) because ties have zero area and the difference quotients are bounded. For a vertex \(X\) strictly inside the hull, the subdifferential of \(-q_X\) is the convex hull of the negative coordinate vectors for the minimizing indices. Thus the minimum condition distributes one unit of weight from \(X\) among the cells with these indices. A vertex in an open cell contributes one unit to that cell. A vertex on an internal dividing line can contribute fractions to several cells. A vertex in the relative interior of a hull side contributes one unit to the cell based on that side. At the minimum, \(\alpha_i\) equals \(c|P\cap Q_i|\) plus \(b\) for the hull side and \(b\) times the total weight from these vertices. Let \(d_i\) count the vertices of \(P\) in \(Q_i\) or in the relative interior of \(E_i\). Each contributes one unit. Vertices on internal dividing lines may contribute additional nonnegative weights, which we discard to obtain
\[
 \alpha_i
        \ge c|P\cap Q_i|+b(1+d_i). \tag{5.4}\label{eq:5.4}
\]
In particular \(\alpha_i>0\), so each cell lies in the triangle on \(E_i\) with apex angle \(\alpha_i\).

\emph{Cuts through the remaining vertices.} Consider a cell \(Q\) with a containing triangle of apex angle \(\theta\) on the same base. Write \(M=|P\cap Q|\) and let \(d\) count the vertices of \(P\) in \(Q\) or in the relative interior of the base. Suppose \(d\ge1\) and
\(\theta\ge b(1+d)+cM\).
Choose a vertex \(Z\) counted by \(d\). A line through \(Z\) in a direction strictly between those of the two legs of the containing triangle meets the relative interior of the base. Let \(0<\tau<\theta\) be the angle measured from the first leg to this direction, in the original Euclidean coordinates. The cut divides the base into two segments. Each is the base of a triangle containing the corresponding piece, with apex angles \(\tau\) and \(\theta-\tau\), respectively. To see this, use affine coordinates in which the containing triangle is \(x,y\ge0,\ x+y\le1\). The cutting direction has two positive components. Its intersections with the lines \(y=0\) and \(x=0\) give the two apices. This works also when \(Z\) lies in the relative interior of the base.

Let \(M_L,M_R\) be the areas inside \(P\) and \(d_L,d_R\) the vertex counts for the two pieces, in order of apex angle \(\tau,\theta-\tau\). The areas vary continuously with \(\tau\). The counts are lower semicontinuous because a previously counted vertex is counted in a piece exactly when it lies strictly on the corresponding side of the cutting line. This includes vertices on the open base. Vertices on old dividing lines stay on the boundary of each piece containing them, and the relative interiors of the new bases lie in that of the old base. Since \(Z\) is on the cut,
\[
 M_L+M_R=M,\qquad d_L+d_R\le d-1.
\]
The sets of directions for which the left and right inequalities fail are
\[
 \{\tau:\tau<b(1+d_L)+cM_L\},\qquad
 \{\tau:\theta-\tau<b(1+d_R)+cM_R\}.
\]
Both sets are open. They are disjoint, since adding the two strict inequalities would give \(\theta<b(2+d_L+d_R)+cM\le b(1+d)+cM\). The first contains an interval near \(0\), and the second contains one near \(\theta\), since \(b>0\). They cannot cover the connected interval \(0<\tau<\theta\). A cut between them therefore satisfies both angle inequalities.

Each cut removes at least one vertex from the count. Starting with \(m\) hull cells, at most \(k-m\) cuts are needed. Discard cells whose intersection with \(P\) has zero area. The at most \(k\) remaining cells have bases on hull sides and lie in triangles on those bases. The triangle containing a cell \(Q\) has apex angle at least \(b+c|P\cap Q|\). The intersections of these cells with \(P\) partition \(P\) up to their boundaries.

\emph{Components inside the polygon.} We use the chord construction of Bogosel, Bucur and Fragalà \cite[Theorem 5 and §6]{ref6}. Fix one of the remaining cells \(Q\), with containing triangle of apex angle \(\theta\). No vertex of \(P\) lies in \(Q\) or in the relative interior of its base. Each side of \(P\) meeting \(Q\) therefore crosses it along a full chord. The open chords are disjoint, and none ends in the open base, since the base lies on a supporting line of the hull and contains no polygon vertex there. The chords divide \(Q\) into convex pieces, each wholly inside or outside \(P\).

If the piece next to the base lies in \(P\), its base lies on \(\partial P\), since points in this piece approach the entire base and \(P\) lies on one side of the base line. For any other piece inside \(P\), join an interior point to a point of the open base, choosing the segment to avoid the chord endpoints. The first chord crossed is a whole side of the piece, since no other chord meets the relative interior of that chord. It separates the whole piece from the whole base and lies on \(\partial P\).

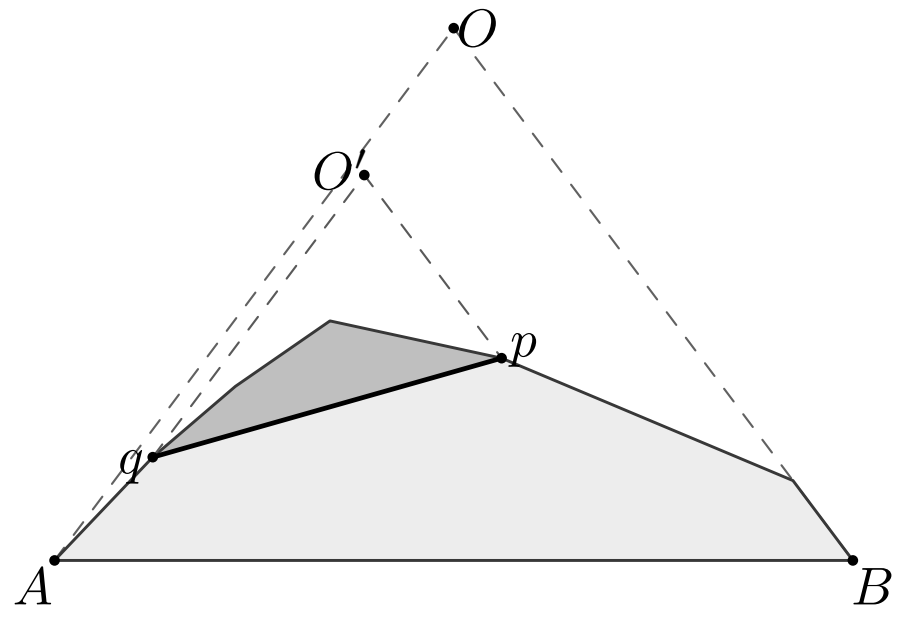
\begin{figure}[htbp]
\centering
\begin{tikzpicture}[x=0.88cm,y=0.88cm,line cap=round,line join=round]
\path[fill=black!7] (3.0000,0.0000) -- (2.5500,0.6000) -- (0.3600,1.5200) -- (-0.9300,1.8000) -- (-1.6410,1.3080) -- (-2.2620,0.7760) -- (-3.0000,0.0000) -- cycle;
\path[fill=black!25] (0.3600,1.5200) -- (-0.9300,1.8000) -- (-1.6410,1.3080) -- (-2.2620,0.7760) -- cycle;
\path[draw=black!62,thin,dashed] (-3.0000,0.0000) -- (0.0000,4.0000) -- (3.0000,0.0000);
\path[draw=black!78,line width=.55pt] (3.0000,0.0000) -- (2.5500,0.6000) -- (0.3600,1.5200) -- (-0.9300,1.8000) -- (-1.6410,1.3080) -- (-2.2620,0.7760) -- (-3.0000,0.0000) -- cycle;
\path[draw=black!65,thin,dashed] (-2.2620,0.7760) -- (-0.6720,2.8960) -- (0.3600,1.5200);
\path[draw=black,line width=.9pt] (0.3600,1.5200) -- (-2.2620,0.7760);
\fill (-3.0000,0.0000) circle (1pt);
\node[font=\small,inner sep=0pt] at (-3.1642,-0.1940) {$A$};
\fill (3.0000,0.0000) circle (1pt);
\node[font=\small,inner sep=0pt] at (3.1493,-0.1940) {$B$};
\fill (0.0000,4.0000) circle (1pt);
\node[font=\small,inner sep=0pt] at (0.1791,4.0000) {$O$};
\fill (0.3600,1.5200) circle (1pt);
\node[font=\small,inner sep=0pt] at (0.5242,1.5797) {$p$};
\fill (-2.2620,0.7760) circle (1pt);
\node[font=\small,inner sep=0pt] at (-2.4262,0.7014) {$q$};
\fill (-0.6720,2.8960) circle (1pt);
\node[font=\small,inner sep=0pt] at (-0.8511,2.9408) {$O'$};
\end{tikzpicture}
\caption{The darker part of the cell lies across \(pq\) from the base \(AB\) and inside the triangle \(O'pq\). The legs through \(O'\) are parallel to the legs through \(O\), so the apex angles are equal.}\label{fig:opening}
\end{figure}

Let \(pq\) be such a chord. Denote the apex and base endpoints of the containing triangle by \(O,A,B\), respectively, and use affine coordinates with \(O=(0,0)\), \(A=(1,0)\) and \(B=(0,1)\), as in Figure~\ref{fig:opening}. The line through \(p,q\) separates some points of \(Q\) from \(AB\). Thus \(O\) is strictly on the side opposite \(AB\), since otherwise the whole triangle would be on the same side as \(AB\). Write the line as \(ux+vy=1\), with \(u,v\ge1\), and label the endpoints so that \(p_x<q_x\) and \(p_y>q_y\). For \(Z\) in the part of \(Q\) across \(pq\) from \(AB\), both \(ZB\) and \(ZA\) meet the chord by convexity. The first intersection gives \(Z_x\ge p_x\) and the second gives \(Z_y\ge q_y\). Together with \(uZ_x+vZ_y\le1\), these inequalities put \(Z\) in the triangle with base \(pq\) and apex \(O'=(p_x,q_y)\). The legs through \(O'\) are parallel to the original legs, so its angle in the original Euclidean coordinates is \(\theta\).

Enumerate the pieces inside \(P\) as \(V_j\). Each piece from \(Q\) lies in a triangle of apex angle \(\theta\) on the chosen base in \(\partial P\). Its base angles are at most those of that triangle, so
\[
 \gamma_j\ge\theta\ge b+c|P\cap Q|\ge b+c|V_j|.
\]
The pieces \(V_j\) partition \(P\) up to their boundaries, as required.\qed

\FloatBarrier
\section{Applications}\label{sec:6}

We first study the dependence on the number of sides and apply the resulting inequalities to polygonal partitions. The final subsection returns to the Fourier comparison and derives bounds for inverse moments of the mixed spectrum.

\subsection{The eigenvalue bound for fewer sides}

Let \(n,h\) be integers with \(1\le h\le n-3\). Theorem~\ref{th:1.1} gives \(c_{n-h}\) as the minimum of \(\Phi\) over polygons with at most \(n-h\) sides. With \(d_n\) defined after \eqref{eq:1.3}, the tangent inequality in that proof gives
\[
 c_{n-h}-c_n>h\,d_n .
\]
For regular polygons of equal area, Dahne, Gómez-Serrano and Pech-Alberich proved that the first eigenvalues and their successive ratios decrease with the number of sides \cite{ref12}. Grinfeld and Strang studied expansions in the reciprocal of the number of sides \cite{ref16}. Berghaus, Georgiev, Monien and Radchenko proved the equal-area expansion and computed its coefficients through order fourteen \cite{ref2}. Here we use Section~\ref{sec:3} to derive the leading coefficient and the derivative estimates for the interpolation \(g(x)=xF(2\pi/x)\), \(x\ge3\).

\Needspace{8\baselineskip}
\begin{corollary}\label{th:6.1}
Let \(J=\lambda_1(\mathbb D)\) as in Proposition~\ref{th:3.1}. As \(n\to\infty\),
\[
 d_n=12\pi J\zeta(3)n^{-4}+O(n^{-6}).
\]
For integers \(1\le h=h(n)\le n-3\), the ratio \((c_{n-h}-c_n)/(h\Delta_n)\) tends to one as \(n\to\infty\) if and only if \(h/n\to0\). For fixed \(h\),
\[
 \frac{c_{n-h}-c_n}{h\Delta_n}
 =1+\frac{2(h-1)}n+\frac{2(h-1)(5h-1)}{3n^2}+O_h(n^{-3}). \tag{6.1}\label{eq:6.1}
\]
\end{corollary}

\noindent\textit{Proof.}\enspace In Section~\ref{sec:3}, \(D(r),H(r),Q(r)\) are bounded as \(r\downarrow0\). Hence the holomorphic extension \(\mathcal L\) of the first mixed eigenvalue in \eqref{eq:3.2} satisfies \(\mathcal L(\sigma)=J+O(r^5)\) uniformly for \(|\sigma-r|\le r/3\). The area formula \eqref{eq:3.1} gives \(Z(\sigma)=1+4\zeta(3)\sigma^3+O(\sigma^5)\). Thus the extension of \(F(2\pi\sigma)\) given by \eqref{eq:3.2} has the form
\[
 F(2\pi\sigma)=\pi J[\sigma+4\zeta(3)\sigma^4]+E(\sigma),
 \qquad \sup_{|\sigma-r|\le r/3}|E(\sigma)|=O(r^6).
\]
Cauchy's estimate gives \(E'(r)=O(r^5)\). Setting \(r=1/x\) gives, as \(x\to\infty\),
\[
 \begin{aligned}
 g(x)&=\pi J[1+4\zeta(3)x^{-3}+O(x^{-5})],\\
 g'(x)&=-12\pi J\zeta(3)x^{-4}+O(x^{-6}).
 \end{aligned}
\]
For integer \(n\ge4\), \(g(n)=c_n\) and \(-g'(n)=d_n\). Define \(B(x)\) by \(-g'(x)=12\pi J\zeta(3)x^{-4}B(x)\). Cauchy's estimate also gives \(E''(r)=O(r^4)\), so \(B(x)=1+O(x^{-2})\) and \(B'(x)=O(x^{-3})\). Put \(m=n-h\). Integrating the bound for \(B'\) gives \(B(x)/B(n)=1+O(m^{-2}-n^{-2})\) uniformly for \(m\le x\le n\) as \(m\to\infty\). With \(v=m/n\), integration of \(-g'\) gives
\[
 \frac{c_m-c_n}{h d_n}
 =\frac{1+v+v^2}{3v^3}\big[1+O(m^{-2}-n^{-2})\big].
\]
For \(0<v\le1\), the factor \((1+v+v^2)/(3v^3)\) tends to one exactly when \(v\to1\). Along any subsequence with \(m\) bounded, \(c_m-c_n\) stays bounded away from zero since \(c_n\) decreases strictly to \(\pi J\), and \(h d_n=O(n^{-3})\). Thus the ratio with \(d_n\) tends to one exactly when \(h/n\to0\). Taking \(h=1\) gives \(\Delta_n/d_n\to1\). For fixed \(h\), expansion at \(v=1\) gives \(1+2h/n+10h^2/(3n^2)+O_h(n^{-3})\) for the ratio with \(d_n\). Dividing this expansion by the one for \(h=1\) gives \eqref{eq:6.1}.\qed

\subsection{Log-convexity for regular polygons}

\begin{theorem}[Log-convexity of the interpolation]\label{th:6.2}
For \(g(x)=xF(2\pi/x)\), the function \(\log g\) is strictly decreasing and has positive second derivative on \([4,\infty)\). The ratios \(c_j/c_{j+1}\) strictly decrease for all integers \(j\ge3\).
\end{theorem}

\noindent\textit{Proof.}\enspace Put \(a=2\pi/x\) and \(d=aF'(a)-F(a)>0\). Since \(g'=-d\) and \(g''=a^2F''(a)/x\), positivity of \((\log g)''\) is equivalent to
\[
                       a^2F(a)F''(a)>d^2.
\]
The bound \(\sin a\,F'(a)\le F(a)\) from the proof of Theorem~\ref{th:1.1} gives
\[
 0<\frac{d}{aF(a)}\le\csc a-\frac1a.
\]
It therefore suffices to prove \(F''(a)/F(a)>(\csc a-1/a)^2\).

For \(0<a\le2\pi/7\), Proposition~\ref{th:3.1} gives \(F(a)\le(8J/15)a\) and \(F''(a)>11J(a/(2\pi))^2/(4\pi)\). Hence
\[
 \frac{F''(a)}{F(a)}>\frac{165a}{128\pi^3}>\frac a{25}
 >\left(\csc a-\frac1a\right)^2.
\]
The second inequality uses \(\pi^3<32\). For the third, \(a<1\) and \(\sin a\ge a-a^3/6\) give \(\csc a-1/a\le a/(6-a^2)<a/5\).

For \(2\pi/7\le a\le\pi/2\), use \(t=\tan(a/2)\), \(z=t^{-2}\), \(G(z)=F(a)/t\), and \(M=2F''(a)/(t(1+t^2))\) as in Section~\ref{sec:4}. The inequality \(a^2F(a)F''(a)>d^2\) becomes
\[
 M>\frac{2d^2}{a^2t^2(1+t^2)G}.
\]
The right side is at most \(2G(z)\cos^2(a/2)(\csc a-1/a)^2\). The identity
\[
 \cot(a/2)\left(\csc a-\frac1a\right)
 =\int_0^1\left(\frac{\sin(va/2)}{\sin(a/2)}\right)^2dv
\]
follows by integrating \(\sin^2(va/2)\). For \(0<v<1\), the integrand increases with \(a\), because \(u\cot u\) decreases on \((0,\pi/2)\). Evaluating at \(a=\pi/2\) and multiplying by \(\sin(a/2)\) gives
\[
 \cos^2(a/2)\left(\csc a-\frac1a\right)^2
 \le\frac{(1-2/\pi)^2}{1+z}.
\]
For \(1\le z\le5/2\), concavity and the values \(G(1)=\pi^2/2\), \(G'(1)=\pi^2/4-2\) give \(G(z)\le G(1)(1+z)/2\). By \eqref{eq:4.5},
\[
 \frac{2G(z)(1-2/\pi)^2}{1+z}\le\frac{(\pi-2)^2}{2}<\frac57<M.
\]
For \(5/2\le z\le9/2\), \eqref{eq:4.1} gives \(G<11/2\). Since \(1-2/\pi<4/11\),
\[
 \frac{2G(z)(1-2/\pi)^2}{1+z}<\frac{11(4/11)^2}{7/2}
 =\frac{32}{77}<\frac{21}{50}<M.
\]
The last bound on \(M\) is again \eqref{eq:4.5}. Since \(\cot^2(\pi/7)<9/2\) by Section~\ref{sec:4}, these inequalities prove \((\log g)''>0\) on \([4,\infty)\).

The equilateral triangle of inradius one has area \(3\sqrt3\) and eigenvalue \(4\pi^2/9\) by Section~\ref{sec:4}, so \(c_3=4\pi^2/\sqrt3\). For the square, \(c_4=2\pi^2\), \(F(\pi/2)=\pi^2/2\), and \(F'(\pi/2)=4\). Thus
\[
 (\log g)'(4)=\frac14-\frac1\pi,\qquad
 \log(2/\sqrt3)\ge1-\sqrt3/2>1/8>1/\pi-1/4.
\]
The chord joining \((3,\log c_3)\) and \((4,\log c_4)\) has smaller slope than the right derivative at four. Use this chord on \([3,4]\) and \(\log g\) on \([4,\infty)\). The resulting interpolation of \(\log c_j\) is decreasing and convex, and its slopes between consecutive integers strictly increase. Equivalently, \(c_j/c_{j+1}\) strictly decreases. This discrete inequality \cite{ref12} also gives \eqref{eq:1.1} once the fixed-side minimum is known, since \(c_j<\sqrt{c_{j-1}c_{j+1}}\le(c_{j-1}+c_{j+1})/2\) for \(j\ge4\).

For the partition estimates below, set
\[
 \beta=\tfrac12\log(c_5/c_7),\qquad
 \varepsilon=\tfrac12\log(c_5c_7/c_6^2).
\]
Then \(\beta>\varepsilon>0\), and summing the increasing successive slopes on either side of six gives
\[
 \log c_j\ge\log c_6-\beta(j-6)+\varepsilon|j-6|
 \qquad(j\ge3). \tag{6.2}\label{eq:6.2}
\]
\qed

\subsection{Periodic convex polygonal tilings}

Hales proved that the regular hexagonal tiling minimizes perimeter per unit area among partitions of the plane into unit-area regions \cite{ref17}. The spectral honeycomb problem asks for the corresponding minimum of a sum of eigenvalues. For a smooth bounded planar domain \(\Omega\) of unit area, Geng and Lin proved that the minimum of \(N^{-2}\sum_{i=1}^N\lambda_1(\Omega_i)\) over partitions into \(N\) cells \(\Omega_i\) tends to a constant independent of \(\Omega\) \cite{ref15}. Caffarelli and Lin conjectured that the constant is \(c_6\), the first eigenvalue of the unit-area regular hexagon \cite{ref8}. In this problem, the boundaries between cells can be curved \cite{ref9,ref20}. For tilings by translates of a single convex polygon, Cesaroni, Fragalà and Novaga proved optimality of the hexagon by Steiner symmetrization \cite{ref10}. We allow finitely many orbits of tiles under the period lattice.

\begin{theorem}[Periodic convex tilings]\label{th:6.3}
Let a locally finite tiling of the plane by bounded convex polygons be invariant under a translation lattice of fundamental area \(A\). Choose one tile \(P_i\) in each of its \(N\) orbits. Write \(j_i\) for the number of sides of \(P_i\) and \(a_i=N|P_i|/A\) for its area divided by the mean tile area. Let \(\varepsilon\) be the constant defined immediately before \eqref{eq:6.2}. Then the normalized geometric mean \(\mathcal G=(A/N)(\prod_i\lambda_1(P_i))^{1/N}\) satisfies
\[
 \log\frac{\mathcal G}{c_6}\ \ge\
 \frac1N\sum_{i=1}^N\left[\varepsilon|j_i-6|
 +\log\frac{\Phi(P_i)}{c_{j_i}}+a_i-1-\log a_i\right]. \tag{6.3}\label{eq:6.3}
\]
Every term on the right is nonnegative. The minimum of \(\mathcal G\) is \(c_6\), attained precisely by regular hexagonal tilings. The normalized power means \((A/N)(N^{-1}\sum_i\lambda_1(P_i)^p)^{1/p}\) for \(p>0\) and the normalized maximum \((A/N)\max_i\lambda_1(P_i)\) have the same minimum and equality case.
\end{theorem}

\noindent\textit{Proof.}\enspace Count vertices and edges on the quotient torus, enlarging the period if necessary so that every tile projects injectively. This multiplies \(A\) and \(N\) by the same factor and preserves the normalized quantities. Split edges at every tile corner. At least three edges meet at each vertex, since two convex faces at a vertex of degree two would both be straight there. For \(V\) vertices and \(E\) edges, Euler's formula and the sum of degrees give \(V-E+N=0\) and \(3V\le2E\). Thus \(\sum_i j_i\le2E\le6N\), including when a tile corner lies in the interior of a neighboring side. Since \(\sum_i a_i=N\),
\[
 \log\frac{\mathcal G}{c_6}
 =\frac1N\sum_i\left[\log\frac{c_{j_i}}{c_6}
 +\log\frac{\Phi(P_i)}{c_{j_i}}+a_i-1-\log a_i\right].
\]
Apply \eqref{eq:6.2} to the first term in the brackets and use \(\sum_i j_i\le6N\) to get \eqref{eq:6.3}. If \(\mathcal G=c_6\), every term on the right of \eqref{eq:6.3} vanishes. By Theorem~\ref{th:1.1}, every tile is a regular hexagon, and \(a_i=1\) for all \(i\). Then \(\sum_i j_i=6N=2E\), so no side has been split in the graph. A regular hexagon determines its neighbor across a whole side, hence the tiling is regular hexagonal. In this tiling all first eigenvalues agree and all the normalized means equal \(c_6\). For any tiling, each positive power mean and the maximum dominate the geometric mean.\qed

The fraction of nonhexagonal tiles is at most \(\log(\mathcal G/c_6)/\varepsilon\). Also
\[
 \frac1N\sum_i|a_i-1|\le 2\sqrt{\log(\mathcal G/c_6)} .
\]
Indeed, for \(a>0\), integrating \(1-1/t\) between \(1\) and \(a\) gives \(a-1-\log a\ge(a-1)^2/(2\max(1,a))\). Apply Cauchy–Schwarz and \(N^{-1}\sum_i\max(1,a_i)\le2\).

For \(\rho>0\), let \(\delta_6(\rho)\) be the infimum of \(\log(\Phi(P)/c_6)\) over unit-area convex hexagons \(P\) at Hausdorff distance at least \(\rho\) from every unit-area regular hexagon. Use \(+\infty\) when the class is empty. To see that \(\delta_6(\rho)>0\), suppose the infimum is zero and take a minimizing sequence. If a hexagon has diameter \(D\), any triangle formed by two vertices at this distance and a third vertex has area at most one. Hence all vertices lie within \(2/D\) of the line through the chosen pair, and the hexagon lies in a strip of width at most \(4/D\). The strip's Dirichlet inequality bounds the diameters in the sequence because their eigenvalues are bounded.

After translation and passage to a subsequence, the vertices converge to those of a convex polygon of area one with at most six sides, allowing coincident vertices. Place the origin inside this limit. The support function of the limit is bounded below by the radius of an interior ball about the origin. Uniform convergence of support functions then puts each member of the subsequence between two dilates of the limit, with dilation factors tending to one. Domain monotonicity and scaling give eigenvalue convergence. By Theorem~\ref{th:1.1}, the limit is regular, contradicting the distance bound. Thus by \eqref{eq:6.3}, the fraction of tiles that have a side count different from six or lie, after rescaling to area one, at Hausdorff distance at least \(\rho\) from every unit-area regular hexagon is at most \(\log(\mathcal G/c_6)/\min\{\varepsilon,\delta_6(\rho)\}\).

\Needspace{10\baselineskip}
\subsection{Finite partitions of a convex polygon}

\begin{theorem}[The hexagonal bound with a boundary term]\label{th:6.4}
Let \(\Omega\) be a convex \(m\)-gon of area \(A\) partitioned into \(N\) positive-area convex polygons \(P_i\) with disjoint interiors and union of closures \(\overline\Omega\). Let \(B\) count the points in the relative interiors of the sides of \(\Omega\) where an interior edge meets the boundary. Use the constants \(\beta,\varepsilon\) defined immediately before \eqref{eq:6.2} and the quantities \(j_i,a_i,\mathcal G\) defined in Theorem~\ref{th:6.3}, now with \(A=|\Omega|\). Then
\[
 \log\frac{\mathcal G}{c_6}\ge\frac{\beta(6-m+B)}{N}
 +\frac1N\sum_{i=1}^N\left[\varepsilon|j_i-6|
 +\log\frac{\Phi(P_i)}{c_{j_i}}+a_i-1-\log a_i\right]. \tag{6.4}\label{eq:6.4}
\]
For fixed \(\Omega\), the infimum of \(\mathcal G\) over these partitions converges to \(c_6\) as \(N\to\infty\). The same limit holds for the infima of the normalized power means of every fixed order \(p>0\) and of the normalized maximum, defined as in Theorem~\ref{th:6.3}.
\end{theorem}

\noindent\textit{Proof.}\enspace Form a plane graph by splitting the cell edges at every cell corner, including at T-junctions, and retaining the \(m\) corners of \(\Omega\). Let \(V,E\) be the numbers of vertices and edges. The graph is connected and its \(N\) bounded faces are the cells, so \(V-E+N=1\). At an interior vertex of degree two, noncollinear edges would leave a reflex angle, while collinear edges would leave both adjacent cells straight. Since the cells are convex and the graph includes only cell corners in the interior, every interior vertex has degree at least three. Each of the \(B\) additional boundary vertices also has degree at least three. The \(m\) corners have degree at least two, so \(2E\ge3V-m\) and \(V\le2N+m-2\). The boundary has \(m+B\) graph edges. Counting interior edges twice and boundary edges once gives
\[
 \sum_i j_i\le 2E-(m+B)\le 6N+m-6-B .
\]
The first inequality allows for sides that were split in the graph. Summing \eqref{eq:6.2} and using the identity for \(\log(\mathcal G/c_6)\) in the proof of Theorem~\ref{th:6.3} gives \eqref{eq:6.4}.

Since \(m\) is fixed and \(B\ge0\), \eqref{eq:6.4} bounds the lower limit of the infima below by \(c_6\). For the upper bound, fix \(0<\delta<1\) and tile the plane by regular hexagons of area \(h=(1-\delta)A/N\). The centers form a triangular lattice. Each hexagon consists of the points for which its center is a nearest lattice point. The hexagons wholly contained in \(\overline\Omega\) cover every point of \(\Omega\) at distance greater than a hexagon's diameter from the boundary. The area of this inner region tends to \(A\) as \(N\to\infty\), so it exceeds \(Nh\) for sufficiently large \(N\).

Choose \(N\) centers of hexagons contained in \(\overline\Omega\) and assign each point of \(\overline\Omega\) to the nearest chosen center, with ties on the boundaries. Each resulting cell is a convex polygon containing the original hexagon at its center. Domain monotonicity gives \(\lambda_1(P_i)\le c_6/h\). In particular the normalized maximum is at most \(c_6/(1-\delta)\). Letting \(\delta\downarrow0\) proves the limits for the maximum, the geometric mean, and the positive power means between them.\qed

For partitions of a fixed \(\Omega\) with \(N\to\infty\) and \(\mathcal G\to c_6\), \eqref{eq:6.4} and the area estimate after \eqref{eq:6.3} show that the proportion of nonhexagons, \(N^{-1}\sum_i|a_i-1|\), and \(B/N\) all tend to zero. For each fixed \(\rho>0\), the fraction of cells that have a side count different from six or lie, after rescaling to area one, at Hausdorff distance at least \(\rho\) from every unit-area regular hexagon is at most
\[
 \frac{\log(\mathcal G/c_6)-\beta(6-m+B)/N}{\min\{\varepsilon,\delta_6(\rho)\}} .
\]

\Needspace{7\baselineskip}
\begin{corollary}[Bends in polygonal partitions]\label{th:6.5}
Let \(P_1,\ldots,P_N\) be positive-area simple polygons with disjoint interiors whose closures have union \(\overline\Omega\), where \(\Omega\) is a convex \(m\)-gon. Use \(A,B,\mathcal G,\beta,\varepsilon\) as in Theorem~\ref{th:6.4}, and let \(j_i\) be the number of sides of \(P_i\). Split the edges at all corners and junctions. Let \(s\) count the interior vertices at which just two edges meet and change direction, and put \(q=(6-m+B-2s)/N\). Then
\[
 \mathcal E_N:=\frac{A}{N^2}\sum_{i=1}^N\lambda_1(P_i)
 \ \ge\ \mathcal G\ \ge\
 c_6\exp\bigl\{\beta q+\varepsilon|q|\bigr\}. \tag{6.5}\label{eq:6.5}
\]
For fixed \(\Omega\), every sequence with \(s=o(N)\) satisfies \(\liminf_{N\to\infty}\mathcal E_N\ge c_6\), and the constant is sharp. If \(\mathcal E_N\le c_6e^{-\eta}\) along a sequence with \(N\to\infty\) and fixed \(\eta>0\), then
\[
 \liminf_{N\to\infty}\frac{s}{N}\ge\frac{\eta}{2\log(c_6/c_7)} .
\]
\end{corollary}

\noindent\textit{Proof.}\enspace Omit vertices of degree two on straight segments. A polygonal path in \(\Omega\) between any two cell interiors can be chosen to avoid vertices and to cross edges transversely. The cells it visits have connected boundaries, and consecutive cells share an edge. Hence the union of the boundaries is connected. With \(V\) vertices and \(E\) edges, the sum of degrees is at least \(3V-m-s\). Euler's formula \(V-E+N=1\) and the \(m+B\) boundary edges give
\[
 \sum_i j_i\le 2E-(m+B)\le 6N+m-6-B+2s .
\]
Put \(D=N^{-1}\sum_i(6-j_i)\). Then \(D\ge q\). The shape and area terms in the identity for \(\log(\mathcal G/c_6)\) are nonnegative by Theorem~\ref{th:1.1} and \(a-1-\log a\ge0\) for \(a>0\). Equation \eqref{eq:6.2} and the triangle inequality give
\[
 \log(\mathcal G/c_6)\ge\beta D+\frac{\varepsilon}{N}\sum_i|6-j_i|
 \ge\beta D+\varepsilon|D|\ge\beta q+\varepsilon|q|.
\]
The last inequality holds because \(\beta>\varepsilon\) makes \(x\mapsto\beta x+\varepsilon|x|\) increasing. The arithmetic mean dominates the geometric mean, proving \eqref{eq:6.5}. For fixed \(m\), \(B\ge0\) gives a lower limit of at least \(c_6\) when \(s=o(N)\). If \(\mathcal E_N\le c_6e^{-\eta}\), then \(q<0\) and \(2s\ge6-m+B+\eta N/(\beta-\varepsilon)\). This gives the stated bound since \(\beta-\varepsilon=\log(c_6/c_7)\). The convex partitions constructed in Theorem~\ref{th:6.4} have \(s=0\) and attain the limiting value.\qed

\begin{remark}
Bending one side of a regular hexagon can lower its first Dirichlet eigenvalue while preserving its area and all six corner angles.

Let \(u\) be the positive ground state on \(R_6\), normalized by \(\int_{R_6}u^2=1\), and put \(\lambda=\lambda_1(R_6)\). On an open side with outward unit normal \(\nu\), the square \((\partial_\nu u)^2\) is nonconstant. If it were constant, Hopf's lemma would make \(\partial_\nu u\) constant there as well. Odd reflection and uniqueness for analytic Cauchy data would force \(u\) to be a constant multiple of \(\sin(\sqrt\lambda\,r)\), where \(r\) is signed distance to the line through the side. By analyticity the formula would hold throughout \(R_6\), contradicting the zero trace on a nonparallel side.

Subtract two smooth bumps of equal integral with compact support in the open side to obtain \(f\) with \(\int f=0\) and \(\int f(\partial_\nu u)^2>0\). A normal displacement by \(tf\) changes the area by exactly \(t\int f=0\) and fixes the side near its endpoints.

Extend \(f\nu\) to a smooth field \(V\) that vanishes near the other sides, and transport \(u\) by \(X\mapsto X+tV(X)\). At \(t=0\), the derivative of the quotient is the integral of \((|\nabla u|^2-\lambda u^2)\operatorname{div}V-2\nabla u\cdot DV\,\nabla u\). This integrand is minus the divergence of the field in Appendix~\ref{sec:A}. Its boundary flux gives derivative \(-\int f(\partial_\nu u)^2<0\). For small \(t>0\), the quotient bounds the new first eigenvalue strictly below \(\lambda\).
\end{remark}

\emph{Convex bodies that need not fill the domain.} For a bounded Lipschitz planar domain \(\Omega\), let \(\mathcal C_k(\Omega)\) be the set of tuples \((K_1,\ldots,K_k)\) of convex bodies contained in \(\Omega\) with pairwise disjoint interiors. We take the Dirichlet eigenvalue on the interior \(K_i^\circ\) of each body. The bodies need not fill \(\Omega\). For every \(p>0\), the theorem of Bucur, Fragalà, Velichkov and Verzini \cite{ref7} now gives
\[
 \lim_{k\to\infty}\frac{|\Omega|^p}{k^{p+1}}
       \inf_{(K_i)\in\mathcal C_k(\Omega)}\sum_{i=1}^k\lambda_1(K_i^\circ)^p=c_6^p,\qquad
 \lim_{k\to\infty}\frac{|\Omega|}{k}
       \inf_{(K_i)\in\mathcal C_k(\Omega)}\max_i\lambda_1(K_i^\circ)=c_6 .
\]
Indeed, \(\lambda_1^p\) decreases under inclusion and scales with degree \(-2p\). By Theorem~\ref{th:1.1} its minimum on unit-area \(j\)-gons is \(c_j^p\). The criterion for the sum requires that the mean of \(c_{j_i}^{p/(p+1)}\) be at least \(c_6^{p/(p+1)}\) for any list of side counts \(j_i\) of mean at most six. This is Jensen's inequality for the decreasing convex interpolation obtained by exponentiating \(p/(p+1)\) times the interpolation of \(\log c_j\) from Theorem~\ref{th:6.2}. For the maximum, the same argument uses \(c_j\). The sharp spectral lower bound for unrestricted partitions remains open.

\subsection{Straight cuts and nodal domains}

\Needspace{8\baselineskip}
\begin{corollary}[Partitions by successive straight cuts]\label{th:6.6}
Suppose a convex \(m\)-gon \(\Omega\) is divided into \(k\) cells \(P_1,\ldots,P_k\) by repeated straight cuts, each joining two boundary points of one current cell. Let \(C\) be the piecewise linear interpolation of \(C(j)=c_j\) for integers \(j\ge3\). Then
\[
 |\Omega|\max_i\lambda_1(P_i)\ge k\,C\!\left(4+\frac{m-4}{k}\right). \tag{6.6}\label{eq:6.6}
\]
The same lower bound holds for \(|\Omega|\lambda\) if a real Dirichlet eigenfunction of eigenvalue \(\lambda\) has \(k\) nodal domains and its entire interior zero set is a finite union of straight segments.
\end{corollary}

\noindent\textit{Proof.}\enspace Each cut adds two copies of the cutting segment and splits at most two existing sides. Thus, with \(j_i\) the number of sides of \(P_i\), we have \(\sum_i j_i\le m+4(k-1)\), including cuts through existing vertices. The slopes of \(C\) are \(-\Delta_{j+1}\), which are negative and strictly increasing by the proof of Theorem~\ref{th:1.1}. The polygonal inequality and Jensen's inequality give
\[
 |\Omega|\max_i\lambda_1(P_i)\ge\sum_i\Phi(P_i)
 \ge\sum_i c_{j_i}\ge k C\!\left(\frac1k\sum_i j_i\right).
\]
Since \(C\) decreases, the side count proves \eqref{eq:6.6}. For a triangle the right side of \eqref{eq:6.6} is \(2\pi^2 k+c_3-c_4\), and for a quadrilateral it is \(2\pi^2 k\). For \(m\ge5\) and \(k\ge m-4\), it is \(2\pi^2 k-(m-4)(c_4-c_5)\).

For the nodal statement, let \(u\) be the eigenfunction. The function \((x,t)\mapsto u(x)\cosh(\sqrt\lambda\,t)\) is harmonic on \(\Omega\times\mathbb R\), so \(u\) is real analytic. Its restriction to a line containing an interior zero segment therefore vanishes along the entire chord. Insert these chords successively, cutting every cell that they cross. The resulting cells are the nodal domains. On each, \(u\) has a fixed sign and zero trace, so \(\lambda\) is its first Dirichlet eigenvalue. Applying \eqref{eq:6.6} proves the nodal bound.\qed

On a fixed convex polygon \(\Omega\), choose a real eigenfunction for each \(\lambda_j(\Omega)\) and let \(k_j\) be its number of nodal domains. Along any subsequence \(j\to\infty\) satisfying the hypothesis on the zero set, Weyl's law \(\lambda_j(\Omega)/j\to4\pi/|\Omega|\) and \eqref{eq:6.6} give \(\limsup k_j/j\le2/\pi\).

For this Weyl asymptotic, use a square grid with spacing \(\delta>0\). The sum of Dirichlet counts on squares contained in \(\Omega\) bounds the Dirichlet count on \(\Omega\) from below. The sum of Neumann counts on squares meeting \(\Omega\) bounds it from above, by zero extension and cutting along grid lines. Counting lattice points in a quarter disk gives \(\delta^2\lambda/(4\pi)+O_\delta(\sqrt\lambda+1)\) eigenvalues below \(\lambda\) on each square for either condition. Divide by \(\lambda\) and let \(\lambda\to\infty\). As \(\delta\downarrow0\), the inner and outer areas both tend to \(|\Omega|\).

On the unit square, \(\sin(\ell\pi x)\sin(\ell\pi y)\) has eigenvalue \(2\pi^2\ell^2\) and \(\ell^2\) nodal domains for each positive integer \(\ell\). It attains \eqref{eq:6.6} for each \(\ell\) and the limiting constant as \(\ell\to\infty\). Polterovich singled out this constant for planar nodal counts \cite[Remark 2.5]{ref21}.

\subsection{Inverse moments of the mixed spectrum}

\Needspace{8\baselineskip}
\begin{corollary}[Mixed spectral zeta values]\label{th:6.7}
Let \(Q\) be a convex polygon with Dirichlet conditions on one side and Neumann conditions on the remaining sides. Suppose the sum of the two interior angles adjacent to the Dirichlet side is \(\pi-\gamma\), where \(0<\gamma<\pi\). For every integer \(r\ge2\),
\[
 \sum_{j\ge1}\frac{1}{\bigl(|Q|\lambda_j^{\mathrm{DN}}(Q)\bigr)^r}
 \le
 \sum_{j\ge1}\frac{1}{\bigl(|I_\gamma|\lambda_j^{\mathrm{DN}}(I_\gamma)\bigr)^r}.
 \tag{6.7}\label{eq:6.7}
\]
Both series converge. Equality for any one \(r\) holds precisely when \(Q\) is similar to \(I_\gamma\) with corresponding Dirichlet side.
\end{corollary}

\noindent\textit{Proof.}\enspace Use the densities \(p,p_*\) from Proposition~\ref{th:2.1}. With the polarized energy and mass forms of Lemma~\ref{th:2.2}, taken linear in the first argument, let \(T_p\) satisfy \(E(T_pu,v)=M_p(u,v)\) on the strip. The change of variables in \eqref{eq:2.2} preserves the energy and the normalized mass on the whole form domain. By the min–max principle the eigenvalues of \(T_p\) are \((|Q|\lambda_j^{\mathrm{DN}}(Q))^{-1}\). The same construction applies to \(p_*\).

To see that both operators are Hilbert–Schmidt, triangulate each polygon and impose Neumann conditions on all the sides of the triangles. This enlarges the form domain, so the number of mixed eigenvalues below \(L>0\) is bounded by the sum of the Neumann counts on the triangles. An invertible affine map bounds each triangle's Rayleigh quotient below by a positive constant times the quotient on a fixed right isosceles triangle. Even reflection across the hypotenuse embeds the latter Neumann problem in a square. Counting the square's cosine modes bounds each count by a constant times \(1+L\). Thus the \(j\)th mixed eigenvalue on either polygon is bounded below by a positive constant times \(j\) for all sufficiently large \(j\). The reciprocal eigenvalues are therefore square summable.

Use the Fourier coefficients \(U_m\) of Lemma~\ref{th:2.2}. Put \(\omega_m(\xi)=\xi^2+(m+\tfrac12)^2\). The energy is \(\pi\sum_m\int_{\mathbb R}\omega_m(\xi)|U_m(\xi)|^2\,d\xi\). Multiplying \(U_m\) by \(\sqrt{\pi\omega_m}\) identifies the energy space with \(L^2\) on \(\{0,1,\ldots\}\times\mathbb R\). For \(X=(m,\xi)\), use \(dX\) for counting measure in \(m\) times Lebesgue measure in \(\xi\). In these coordinates, for \(Y=(n,\eta)\) the kernel of \(T_p\) is
\[
 K_p(X,Y)=
 \frac{\widehat p(\eta-\xi,m-n)+\widehat p(\eta-\xi,m+n+1)}
      {2\pi^2\sqrt{\omega_m(\xi)\omega_n(\eta)}} .
\]
On a set of finitely many modes and a bounded frequency interval the displayed kernel is bounded. The cosine product formula gives its integral against every product of two test functions supported on that set. Finite sums of such products are dense in \(L^2\) of the pair of variables, so the formula agrees almost everywhere there with the Hilbert–Schmidt kernel. Exhaustion gives the formula on the whole space. Let \(K_*\) be the kernel with \(p_*\) in place of \(p\). Proposition~\ref{th:2.1} and \(\widehat p_*>0\) give \(|K_p(X,Y)|\le K_*(X,Y)\). Unless \(Q\) is similar to \(I_\gamma\) with corresponding base, the bound is strict for every \(X,Y\) because the frequency with angular index \(m+n+1\) is nonzero.

\Needspace{5\baselineskip}
For each integer \(R\ge1\), set \(\mathcal X_R=\{0,\ldots,R\}\times[-R,R]\) and let \(P_R\) be multiplication by its indicator. The kernels are bounded on \(\mathcal X_R^2\). Composing and taking the trace, with \(X_{r+1}=X_1\), gives
\[
 \operatorname{tr}\bigl((P_RT_pP_R)^r\bigr)
 =\int_{\mathcal X_R^r}\prod_{i=1}^r K_p(X_i,X_{i+1})\,dX_1\cdots dX_r.
\]
For \(r=2\) the trace is the square of the Hilbert–Schmidt norm. For \(r>2\), use the trace formula for the product of the \((r-1)\)st power and the last factor. The product of the \(K_*\) bounds the real part of the integrand above. When the kernel inequality is strict, the difference after taking real parts has positive integral already on \((\{0\}\times[-1,1])^r\). This lower bound on the difference is independent of \(R\).

As \(R\to\infty\), the truncated kernels converge in \(L^2\), so the operators converge in Hilbert–Schmidt norm. Their \(r\)th powers converge in trace norm. Indeed, expand the difference of powers as a telescoping sum. In each term bound the difference factor and one other factor in Hilbert–Schmidt norm and the rest in operator norm. The trace inequality therefore passes to the limit, with a strictly positive difference in the strict case. This proves \eqref{eq:6.7}. Similarity with the bases corresponding preserves all the area-normalized eigenvalues and gives equality.\qed

For simply connected planar domains of fixed conformal radius, Laugesen and Morpurgo proved that the disk minimizes the Dirichlet spectral zeta function at each real exponent greater than one \cite{ref30}. For a doubly connected membrane, Laugesen averaged the two assignments of Dirichlet and Neumann conditions to the boundary components and bounded the average below by the value for a cylinder of prescribed modulus and mass \cite[Theorem 6]{ref33}. For a weighted cylinder he obtained a lower bound with one assignment under a cumulative mass condition \cite[Theorem 4]{ref34}. Corollary~\ref{th:6.7} gives an upper bound for a convex polygon with one Dirichlet side at prescribed area and sum of base angles.

To recover the first mixed eigenvalue inequality, write \(\mu_j=(|Q|\lambda_j^{\mathrm{DN}}(Q))^{-1}\) in nonincreasing order. For integers \(r\ge2\),
\[
 \mu_1^r\le\sum_j\mu_j^r\le\mu_1^{r-2}\sum_j\mu_j^2 .
\]
The same bounds hold for \(I_\gamma\). Take \(r\)th roots in \eqref{eq:6.7} and let \(r\to\infty\) through integers to obtain \eqref{eq:2.1}. The eigenvalue growth estimate ensures convergence of both series for every real \(r>1\). Does the comparison also hold for noninteger \(r>1\)?

\section{Concluding remarks}\label{sec:7}

At unit area, Brasco, De Philippis and Velichkov bound the excess over the first eigenvalue of the disk below by a constant times the squared area of the symmetric difference from a closest unit-area disk \cite{ref28}.

A Dirichlet test function on a polygon restricts to an admissible function for the mixed problem on each cell. This restriction need not preserve orthogonality to the ground state. For unit-area \(n\)-gons \(P\) and \(R_n\), with \(R_n\) regular, what is the sharp bound on \(\lvert\lambda_j(P)-\lambda_j(R_n)\rvert\) in terms of \(\lambda_1(P)-\lambda_1(R_n)\) for fixed \(j\ge2\)? For the ball, Bucur, Lamboley, Nahon and Prunier obtained sharp bounds for individual eigenvalues and sums over clusters \cite{ref31}.

For partitions of a fixed convex polygon into \(N\) polygonal cells, suppose the normalized sum \(\mathcal E_N\) from \eqref{eq:6.5} stays a fixed amount below \(c_6\) as \(N\to\infty\). Corollary~\ref{th:6.5} then gives at least a positive multiple of \(N\) interior vertices where two edges meet and change direction. The unrestricted problem of Caffarelli and Lin allows curved interfaces \cite{ref8}. In an optimal planar partition, three \(C^{1,\alpha}\) arcs meet at \(120^\circ\) at a junction of Almgren frequency \(3/2\), by the theorem of Ognibene and Velichkov \cite{ref32}. The example following Corollary~\ref{th:6.5} motivates studying how a bend of a shared side changes the sum of the two adjacent eigenvalues.

\appendix
\section{The triangular case}\label{sec:A}

We use the boundary moment and reflection argument of Fragalà and Velichkov \cite[Theorem 1.1]{ref29}.

\begin{proposition}[The minimizing triangle]\label{th:A.1}
Among triangles of prescribed positive area, the equilateral triangle uniquely minimizes the first Dirichlet eigenvalue up to congruence.
\end{proposition}

\noindent\textit{Proof.}\enspace A triangle of area \(A\) whose longest side has length \(D\) lies in a strip of width \(2A/D\). The Dirichlet inequality across the strip gives \(\lambda_1\ge\pi^2D^2/(4A^2)\). Thus a minimizing sequence of fixed area has bounded diameter. After rigid motions and passage to a subsequence, its vertices converge to those of a triangle \(T\) of the same area. The affine maps from \(T\) to the triangles in the subsequence converge to the identity. Comparing Rayleigh quotients under these maps proves that \(T\) is a minimizer.

Choose any side as the base and write the vertices as \((-b,0),(b,0),(t,h)\), with \(b,h>0\). Let \(u>0\) be the ground state, normalized by \(\int_Tu^2=1\), and put \(\lambda=\lambda_1(T)\). Write \(\nu\) for the outward unit normal on each open side. A variation that tilts the base about its midpoint gives, along a sequence of radii \(\rho\downarrow0\),
\[
 \lim_{\rho\downarrow0}\int_{-b+\rho}^{b-\rho}x\,u_y(x,0)^2\,dx=0 . \tag{A.1}\label{eq:A.1}
\]
Take the affine field \(V\) that vanishes at the apex and moves the two base endpoints along their incident legs with vertical velocities \(b\) on the left and \(-b\) on the right. The field is tangent to both legs and satisfies \(V\cdot\nu=x\) on the base. Its divergence is constant with integral zero, so \(\operatorname{div}V=0\). For small \(s\), transport \(u\) by \(X\mapsto X+sV(X)\). The transported Rayleigh quotient times the area is minimized at \(s=0\) because \(\Phi\) is scale invariant and \(T\) minimizes it. The derivatives of the area and transported \(L^2\) mass vanish since \(\operatorname{div}V=0\). Differentiating the energy gives \(\int_T\nabla u\cdot DV\,\nabla u=0\). The vector field
\[
 2(V\cdot\nabla u)\nabla u-|\nabla u|^2V+\lambda u^2V
\]
has divergence \(2\nabla u\cdot DV\,\nabla u\) by the eigenvalue equation and \(\operatorname{div}V=0\). On an open side its normal component is \((V\cdot\nu)|\partial_\nu u|^2\).

To control the flux near a vertex, let \(B_r\) be the disk of radius \(r\) centered there and \(\ell\) the coordinate along the unit inward bisector, measured from the vertex. The corner angle is less than \(\pi\), so \(\ell\) is bounded above and below by positive multiples of the distance to the vertex. In a small \(T\cap B_r\), a positive multiple of \(\sin(\sqrt\lambda\,\ell)\) dominates \(u\) on the boundary and solves the same equation. The first Dirichlet eigenvalue of \(T\cap B_r\) is strictly greater than \(\lambda\), since equality would give by zero extension a ground state on \(T\) vanishing on an open set. Testing the difference with its negative part proves that \(u\) is bounded by a constant times the distance to the vertex.

Testing with \(u\) times the square of a cutoff supported in \(B_{2r}\) and equal to one on \(B_r\) gives
\[
 \int_{T\cap B_r}|\nabla u|^2
 \le C(r^{-2}+\lambda)\int_{T\cap B_{2r}}u^2=O(r^2),
\]
where \(C\) is independent of \(r\). Average the sum of the integrals of \(|\nabla u|^2\) over the three circular arcs of common radius between \(r\) and \(2r\). Along a common sequence \(\rho\downarrow0\), the sum on the arcs of radius \(\rho\) is \(O(\rho)\), and the flux through these arcs tends to zero. Integrating the divergence over the triangle with the three corner disks removed proves \eqref{eq:A.1}.

If \(t>0\), the midpoint of the horizontal slice at height \(y\) has horizontal coordinate \(ty/h>0\) for \(0<y<h\). Reflecting across \(x=0\) therefore maps \(T\cap\{x<0\}\) into \(T\), sending its open oblique side into the interior. On \(T\cap\{x<0\}\), the difference \(w(x,y)=u(-x,y)-u(x,y)\) solves \(-\Delta w=\lambda w\), with nonnegative boundary values and strictly positive values on the open oblique side. Strict domain monotonicity and the negative-part test give \(w\ge0\), and the strong maximum principle gives \(w>0\). Hopf's lemma on the base, applied to \(w\) and \(u\), yields \(u_y(x,0)>u_y(-x,0)>0\) for \(0<x<b\). The moment in \eqref{eq:A.1} is then bounded below by the integral of \(x[u_y(x,0)^2-u_y(-x,0)^2]\) over any fixed closed subinterval of \((0,b)\) of positive length, for all sufficiently small \(\rho\). This integral is positive, contradicting \eqref{eq:A.1}. Reflection also excludes \(t<0\). Thus the opposite vertex lies on the perpendicular bisector of each side, and every minimizer is equilateral.\qed

\section{A convolution identity for two random points}\label{sec:B}

Let \(Q\) be as in Proposition~\ref{th:2.1}, with base angles \(\pi a,\pi b\) and \(q=1-a-b\). Use the densities \(p,p_*\) on \(\mathcal C\) from that proposition and the measure \(d\mu\) of Section~\ref{sec:2}. Choose two independent uniform points in \(Q\). For each point with strip coordinates \((x,t)\), replace \(t\) by \(t\) or \(-t\) with equal probability, independently for the two points. Each resulting point on \(\mathcal C\) has density \(p\). For any density \(f\) on \(\mathcal C\), write \(\check f(x,t)=f(-x,-t)\). The difference of the two points, taken modulo \(2\pi\) in the angular coordinate, has density \(p*\check p\).

For \(v>0\), put \(h(v)=(e^{-av/2}-e^{-bv/2})^2(1-e^{-qv})/[v(1-e^{-v})]\). This is \(O(v)\) at zero and decays exponentially at infinity. Let \(K_v\) be the law of a horizontal shift by \(v/2\) or \(-v/2\) with equal probability and an independent angular shift. Its angular density is the Poisson kernel \((1-\rho^2)/(1-2\rho\cos\theta+\rho^2)\) for \(d\theta/(2\pi)\), where \(\rho=e^{-v/2}\). The Fourier series gives \(\widehat K_v(\xi,k)=e^{-|k|v/2}\cos(\xi v/2)\).

Let \(M=\int_0^\infty 2h(v)K_v\,dv\). Write \(M^{*j}\) for its \(j\)-fold convolution, with \(M^{*0}\) the point mass at the origin. Then \(\Pi=e^{-M(\mathcal C)}\sum_{j\ge0}M^{*j}/j!\) is a probability measure with transform \(\exp[-2\int_0^\infty h(v)(1-\widehat K_v)\,dv]\). Let \(\nu\) be the probability measure in \eqref{eq:2.3} and \(\check\nu\) its image under \((x,t)\mapsto(-x,-t)\). Equation \eqref{eq:2.6} gives \(|\widehat p_{a,b}|^2=\widehat{p_* *\check p_*}\,\widehat\Pi\). By \eqref{eq:2.3} and uniqueness of Fourier transforms,
\[
             p*\check p=(p_* *\check p_*)*(\nu*\check\nu)*\Pi . \tag{B.1}\label{eq:B.1}
\]
Thus the difference has the same law as the corresponding difference for the isosceles triangle, with independent shifts of laws \(\nu*\check\nu\) and \(\Pi\) added. The eigenvalue comparison in Section~\ref{sec:2} uses the Fourier bound for \(p\) itself.

\microtypesetup{expansion=false}
\section*{Use of AI tools}

OpenAI’s ChatGPT was used to develop parts of the Fourier comparison and the polygonal partition argument, check intermediate estimates and equality cases, search the literature and revise the text. The author takes full responsibility for the mathematical results and the text of this paper.

\end{document}